\documentclass[11pt,a4paper,twoside]{article}
\usepackage{latexsym,amsfonts,amsmath, amsthm, amssymb}

\usepackage{
epsfig, graphicx, psfrag, subfigure }

\newtheorem{theorem}{Theorem}
\newtheorem{proposition}[theorem]{Proposition}
\newtheorem{lemma}[theorem]{Lemma}
\newtheorem{corollary}[theorem]{Corollary}

\theoremstyle{definition}
\newtheorem{example}[theorem]{Example}
\newtheorem{definition}[theorem]{Definition}
\newtheorem{remark}[theorem]{Remark}

\newcommand{\dR}{\ensuremath{\mathbb R}}
\newcommand{\dZ}{\ensuremath{\mathbb Z}}

\newcommand{\var}{\mathrm{Var}} 
\newcommand{\ent}{\mathrm{Ent}} 
\newcommand{\C}[1]{\ensuremath{{\mathcal C}^{#1}}} 
\renewcommand{\L}{\ensuremath{\mathbb L}} 
\newcommand{\GI}{\mathbf{L}} 
\newcommand{\PT}[1]{\mathbf{P}_{\!#1}} 

\newcommand{\NRM}[1]{\ensuremath{{\left\Vert #1\right\Vert}}} 
\newcommand{\NRME}[2]{\ensuremath{{\left\Vert #1\right\Vert_{#2}}}} 

\newcommand{\ind}{\mathrm{1}\hskip -3.2pt \mathrm{I}} 

\begin{document}

\newcommand{\met}{\qquad \mbox{et} \qquad}
\title{Orlicz-Sobolev inequalities for sub-Gaussian measures and ergodicity of Markov semi-groups
\thanks{Supported by EPSRC GR/S61690/01 \& GR/R90994/01} }
\author{C.~Roberto and B.~Zegarli{\'n}ski
\smallskip
\\ {\small Universit{\'e} de Marne la Vall{\'e}e}
\\ {\small Imperial College London}
}
\date{}

\maketitle \pagestyle{myheadings}
\markboth{C. Roberto and B. Zegarli\'nski}
{Orlicz-Sobolev inequalities and ergodicity of Markov semi-groups}

\begin{abstract}
We study coercive inequalities in Orlicz spaces associated
to the probability measures on finite and infinite dimensional spaces
which tails decay slower than the Gaussian ones.
We provide necessary and sufficient criteria for such inequalities to
hold and discuss relations between various classes of inequalities.
\end{abstract}


\section{Introduction}

Sobolev type inequalities play an essential role in the study of the decay
to equilibrium of Markov semi-groups to their associated probability measure.
Several surveys deal with the celebrated Poincar\'e inequality
and the stronger logarithmic Sobolev inequality, see e.g.
\cite{Gross}, \cite{Bak94}, \cite{Led99}, \cite{ane}, \cite{GZ03} and
\cite{Ro99}. It appears that the Poincar\'e inequality is particularly
adapted to the study of the two sided exponential measure while the logarithmic Sobolev inequality is the perfect tool to deal with the Gaussian measure. Both are now well understood.

In recent years intermediate measures, as for example
$$
d \mu_\alpha (x) = (Z_\alpha)^{-1} e^{-|x|^\alpha}
\qquad
\alpha \in (1,2) ,
$$
attracted a lot of attention (note that
for such measures the Logarithmic Sobolev inequalities cannot hold).
To deal with such measures, several authors
generalized the Poincar\'e and logarithmic Sobolev inequalities
in the following way.

Recall first that a probability measure $\mu$, say on $\dR^n$,
is said to satisfy a Poincar\'e inequality
if there exists a constant $C$
such that every $f :\dR^n \to \dR$ smooth enough satisfies
$$
\var_\mu(f) \leq C \int |\nabla f|^2 d\mu
$$
and to satisfy a logarithmic Sobolev inequality if
$$
\ent_\mu(f^2) \leq C \int |\nabla f|^2 d\mu
$$
where $\var_\mu(f)=\mu(f^2) - \mu(f)^2$
is the variance (for short
$\mu(f)=\int f d\mu$), and
$\ent_\mu(f)=\mu \left( f \log (f/\mu(f)) \right)$ is the entropy of a
positive function.

The latter can be rewritten in the form
$$
\int f^2 \log (f^2) d\mu - \int f^2 d\mu \log \left( \int f^2 d\mu \right)
\leq C \int |\nabla f|^2 d\mu
$$
or equivalently
$$
\lim_{p \to 2^-} \frac{ \int f^2 d\mu - \left( \int |f|^p \right)^\frac 2p}{2-p}
\leq 2 C \int |\nabla f|^2 d\mu .
$$
Hence two natural generalizations are the following additive $\Phi$-Sobolev
inequality
\begin{equation} \label{eq:dap} \tag*{(\textbf{$\Phi$-S})}
\int \Phi(f^2) d\mu - \Phi \left( \int f^2 d\mu \right)
\leq C \int |\nabla f|^2 d\mu
\end{equation}
and the  Beckner-type inequality:
\begin{equation} \label{eq:bti}
\sup_{p \in [1,2)} \frac{ \int f^2 d\mu - \left( \int |f|^p \right)^\frac 2p}{T(2-p)}
\leq 2 C \int |\nabla f|^2 d\mu .
\end{equation}
Inequality (\textbf{$\Phi$-S}) has been introduced in
\cite{barthe-cattiaux-roberto} as an intermediate tool to prove an
isoperimetric inequality for the measure $\mu_\alpha$. It is also
related to the work by Chafa\"i \cite{chafai}. On the other hand,
Beckner introduced in \cite{Bec89}  Inequality \eqref{eq:bti} with
$T(r)=r$ in his study of the Gaussian measure $\mu_2$.
Lata{\l}a and Oleszkiewicz \cite{latala-Oleszkiewicz} consider the more general
$T_\alpha(r)=r^{2(1-\frac{1}{\alpha})}$, $\alpha \in (1,2)$ and prove that
$\mu_\alpha$ satisfies Inequality \eqref{eq:bti} with such $T_\alpha$. Furthermore this
inequality appears to be well adapted to the study of concentration
of measure phenomenon via the celebrated Herbst argument. Further
generalizations are done in this direction in
\cite{barthe-cattiaux-roberto}, see also \cite{wang}. When
$T=T_\alpha$, Inequality \eqref{eq:bti} is known as
the Lata{\l}a and Oleszkiewicz
Inequality.

\medskip

While
the logarithmic Sobolev inequality enjoys a lot of properties and applications (tensorisation,
concentration of measure, isoperimetry,
decay to equilibrium, hypercontractivity), none of its generalizations
appears to be well adapted simultaneously to all these properties and applications.
This is the main reason why one has to generalize in different ways the
logarithmic Sobolev inequality.

\medskip

Motivated by this, in this paper we study the following
new generalization we shall call the Orlicz-Sobolev inequality
\begin{equation} \label{in:os}
\tag*{({\textbf{O-S}})}
\NRME{(f-\mu(f))^2}{\Phi} \leq C \int |\nabla f|^2 d \mu
\end{equation}
where the constant $C$ is independent of the function $f$. Here,
$\NRME{\cdot}{\Phi}$ denotes the Luxembourg norm associated
to the Orlicz function $\Phi$ and the probability measure $\mu$
on finite or infinite products of real lines $\dR$.

If $\Phi(x)=|x|^p$, $p \in [1, \infty)$ then
$\NRME{(f-\mu(f))^2}{\Phi}=\NRME{(f-\mu(f))^2}{p}$.
Thus for $1< p <\infty$, inequality ({\textbf{O-S}})
can be considered as a Sobolev type inequality.
For $p=1$ it is the Poincar{\'e} inequality.
On the other hand, for $\Phi(x)=|x| \log (1+|x|)$, it is proved in
\cite{bobkov-gotze} that ({\textbf{O-S}})
is equivalent
(up to universal constants) to the logarithmic Sobolev inequality.
Thus, for an interpolation family of Orlicz functions going from $|x|$ to
$|x| \log (1+|x|)$ (as for instance
$|x| \log(1+|x|)^\beta$, $\beta \in [0,1]$),
({\textbf{O-S}})
is an interpolating family of functional inequalities between
Poincar{\'e} and the logarithmic Sobolev inequality.

\medskip

Our first objective is to give in Section 2 a constructive criterium for a
probability measure on a finite dimensional Euclidean space to satisfy such
an inequality. In particular we will prove (Corollary \ref{cor:alpha})
that the sub-Gaussian probability measures $\mu_\alpha$ (and product of it)
satisfy the Orlicz-Sobolev Inequality ({\textbf{O-S}}) with
$\Phi(x) = |x| \log(1+|x|)^{2(1-\frac{1}{\alpha})}$.

Note that the Orlicz-Sobolev Inequality ({\textbf{O-S}}) need not tensorise
in general. Hence, in order to get dimension free results,
we will use our criterium and ideas from
\cite{barthe-cattiaux-roberto} to prove the equivalence between
the Orlicz-Sobolev Inequality  and the Beckner type Inequality
\eqref{eq:bti} that do tensorise.

Finally, using our results, we prove that under suitable mixing conditions
the Lata{\l}a-Oleszkiewicz inequalities are satisfied for
Gibbs measures on infinite dimensional spaces.
This provides an extension of a result discussed in
\cite{GZ03} to a comprehensive family of local specifications.

\medskip

In Section 3 we discuss the implications of Orlicz-Sobolev inequalities for the
decay to equilibrium in Orlicz norms for Markov semi-group with the generator
given by the corresponding Dirichlet form. This includes in particular a necessary
and sufficient condition for the exponential decay, which extends a well known
classical property of the $\L_2$ space and Poincar{\'e} inequality.
One of our main result states that the Orlicz-Sobolev inequality
(\textbf{O-S}) implies, under mild assumptions on $\Phi$, that
\begin{equation} \label{eq:12}
\NRM{\PT{t}f}_\Phi \leq e^{-ct}  \NRM{f}_\Phi
\end{equation}
for any $f$ with $\mu(f)=0$.
Our technical development allows us to consider at the end of the section
the case of decay to equilibrium for functionals which do not have convexity
property of the norm as for example functionals of the form
$\mu \left(|f|^q\log |f|^q/\mu(|f|^q) \right)$
with $q>1$. In case of relative entropy corresponding to $q=1$
and a hypercontractive diffusion semi-group the exponential decay is well known.
For $q>1$ we show that after certain characteristic period of time
one gets (essentially) exponential decay and by suitable averaging
one can redefine the functional so it has the exponential decay property.

\medskip

In Section 4 we discuss a relation between Orlicz-Sobolev
and the additive $\Phi$-Sobolev (\textbf{$\Phi$-S}) inequalities.
The additive $\Phi$-Sobolev inequalities naturally tensorise.
We show that it also has an analog of the mild perturbation property which allows
to construct local specifications satisfying
such the inequality.
Moreover we prove that, if the local specification is mixing, similar arguments
to those employed in the proof of Logarithmic
Sobolev inequalities work in the current situation.
By this we get a constructive way to provide examples of nontrivial
Gibbs measures on infinite dimensional spaces satisfying
the additive $\Phi$-Sobolev inequalities.
In a forthcoming paper \cite{fougeres-r-z} we will use them in the study of
infinite dimensional nonlinear Cauchy problems.

\medskip

In order to show a decay to equilibrium in a stronger than $\L_2$ sense,
in Section 5 we introduce and study certain natural generalization of
Nash inequalities which follow from Orlicz-Sobolev inequalities.
Such inequalities provide a bound on a covariance in terms of the
Dirichlet form and suitable (weaker than $\L_2$) Orlicz norm.
One illustration of our results is that the Inequality (\textbf{O-S}) proved in section
1 for $\mu_\alpha$ and $\Phi(x)=|x|\log(1+x)^{2(1-\frac{1}{\alpha})}$
implies that the associated semi-group $(\PT{t})_{t \geq 0}$
is a continuous map from $\L_\Psi$ into $\L_2$ with
$\Psi(x)=x^2/\log(1+|x|)^{2(1-\frac{1}{\alpha})}$. Furthermore
$$
\NRM{\PT{t} }_{\L_\psi \to \L_2} \leq \frac{C_{\alpha}}{t^\gamma} \qquad
\qquad \forall t>0
$$
for some positive constant $C_\alpha$ and $\gamma$. This result state that as soon as $t$ is positive,
the semi-group regularizes any initial data from $\L_\Psi$ into $\L_2$.
(For general discussion about the interest and application of Nash-type inequalities,
we refer the reader to e.g.
\cite{davies}, \cite{saloff}, \cite{varopoulos}, \cite{GZ03}.)
Note that this bound is different from \eqref{eq:12} where on both sides appear
the same $\L_\Phi$ norm.

\medskip

As a summary, all the multitude of the inequalities and relations between them discussed in this work is illustrated with the corresponding implication network diagram provided at the end of the paper. Since in our investigations we have used intensively numerous properties of Young functions and Orlicz/Luxemburg norms,  for the convenience of the reader in the Appendix we gathered a plentitude of useful facts.

For other directions on the study of sub-Gaussian measures,
the reader could like to see also
\cite{gentil-guillin-miclo,barthe-cattiaux-roberto2,
barthe-cattiaux-roberto3,BZ01}.


\paragraph{Acknowledgements:}
\textit{The authors would like to thank Pierre Foug{\`e}res
for the critical reading of the manuscript and the referee for his valuable observations.
Cyril Roberto also warmly
acknowledges the hospitality at the
Imperial College.}


\section{A criterium for Orlicz-Sobolev inequalities}

In this section we provide a criterium for inequality ({\textbf{O-S}})
to hold.
This criterium allows us to prove that Orlicz-Sobolev inequalities are
equivalent, up to universal constants, to Bekner-type inequalities. In turn,
we give a family of Orlicz functions for which the Orlicz-Sobolev
inequality holds for a corresponding sub-Gaussian measure.
We end with an application to Gibbs measure on infinite state space.

In \cite{barthe-cattiaux-roberto},
the authors introduce a general tool to
obtain a criterium which is based on an appropriate notion of capacity
(\cite{mazya}) initially
introduced in \cite{barthe-roberto}.
More precisely, let $\mu$ and $\nu$ be two absolutely
continuous measures on $\dR^n$.
Then, for any Borel set $A \subset \Omega$, we set
$$
\mathrm{Cap}_\nu(A,\Omega) :=
\inf \left\{ \int |\nabla f|^2 d\nu ; f \geq \ind_A
\mbox{ and } f_{|\Omega^c}=0 \right\} .
$$
If $\mu$ is a probability measure on $\dR^n$, then,
for $A \subset \dR^n$ such that $\mu(A) < \frac{1}{2}$,
the capacity of $A$ with respect to $\mu$ and $\nu$ is
\begin{eqnarray*}
\mathrm{Cap}_\nu(A,\mu) & := &
\inf \left\{ \int |\nabla f|^2 d\nu ; \ind \geq f \geq \ind_A
\mbox{ and } \mu(f=0) \geq \frac 12 \right\} \\
& = &
\inf \left\{ \mathrm{Cap}_\nu(A,\Omega) ; \Omega \subset \dR^n \text{
    s.t. }
\Omega \supset A
\mbox{ and } \mu(\Omega)= \frac 12 \right\} .
\end{eqnarray*}
For simplicity we will write $\mathrm{Cap}_\mu(A)$ for $\mathrm{Cap}_\mu(A,\mu)$.
[For a general introduction and discussion on the notion of capacity we refer the reader to \cite[section 5.2]{barthe-cattiaux-roberto} .]
The second equality in the above definition comes from the fact that
$\mathrm{Cap}_\nu(A,\Omega)$ is non-increasing in $\Omega$
and a suitable truncation argument
(see \cite{barthe-cattiaux-roberto}).

We start with the following criterium in dimension $n$ and its more explicit form
in dimension one.

\begin{theorem}\label{th:cos}
Let $\mu$ and $\nu(dx)=\rho_\nu(x)dx$ be two absolutely continuous probability measures
on $\dR^n$. Consider a Young function
$\Phi$ and  fix $k \in (0,\infty)$
such that for any function $f$ with $f^2 \in \L_\Phi (\mu)$,
one has $\NRME{\mu(f)^2}{\Phi} \leq k \NRME{f^2}{\Phi}$.
Let $C_\Phi$ be the optimal constant such that for any smooth
function $f : \dR^n \rightarrow \dR$ one has
$$
\NRME{(f-\mu(f))^2}{\Phi}
\leq
C_\Phi \int |\nabla f|^2 d \nu .
$$
Then $\frac{1}{8} B(\Phi) \leq C_\Phi
\leq 8(1+k) B(\Phi)$
where $B(\Phi)$ is the smallest
constant such that for every
$A \subset \dR^n$ with $\mu(A) < \frac12$,
$$
\NRME{\ind_A}{\Phi}  \leq B(\Phi)
\mathrm{Cap}_\nu(A,\mu) .
$$
Moreover if $n=1$, one has
$$
\frac{1}{8} \max(B_+(\Phi),B_-(\Phi))
\ \leq \  C_\Phi \ \leq \
8(1+k) \max(B_+(\Phi),B_-(\Phi))
$$
where
$$
B_+(\Phi) = \sup_{x >m} \NRME{\ind_{[x,+\infty)}}{\Phi}
\int_m^x \frac 1{\rho_\nu} ,
$$
$$
B_-(\Phi) = \sup_{x < m} \NRME{\ind_{(-\infty,x]}}{\Phi}
\int_x^m \frac 1{\rho_\nu} ,
$$
and $m$ is a median of $\mu$.
\end{theorem}

\begin{remark} \label{rem:thcos}
Note that
by the property \eqref{equiv} in the Appendix,
$\NRME{\ind_A}{\Phi} = 1/\Phi^{-1}(1/\mu(A))$. In particular for
$\mu(A)<\frac12$ we have $\NRME{\ind_A}{\Phi}< 1/\Phi^{-1}(2) $.

For explanation concerning the condition
$\NRME{\mu(f)^2}{\Phi} \leq k \NRME{f^2}{\Phi}$
when $f^2 \in \L_\Phi (\mu)$, see Lemma \ref{lemma_injection_L1} and
Remark \ref{remark_norm_mu_f2} in the Appendix.
\end{remark}

\begin{proof}
Fix a locally Lipschitz function
$f : \dR^n \rightarrow \dR$ and let $c$ be a median of $f$, {\it i.e.}
$\mu(f \geq c) \geq \frac 12$ and $\mu(f \leq c) \geq \frac 12$.
Then define $f_+=(f-c) \ind_{f > c}$ and $f_-=(f-c) \ind_{f < c}$.
By assumption about $\Phi$,
\begin{eqnarray*}
\NRME{(f-\mu(f))^2}{\Phi}
&=&
\NRME{(f-c + \mu(f-c))^2}{\Phi} \\
&\leq &  2 \NRME{(f-c)^2}{\Phi} + 2 \NRME{\mu(f-c)^2}{\Phi} \\
&\leq &
2(1+k) \NRME{(f-c)^2}{\Phi} \\
&\leq & 2(1+k)(\NRME{f_+^2}{\Phi}+\NRME{f_-^2}{\Phi}) .
\end{eqnarray*}
with $k\in(0,\infty)$ independent of $f$.
It follows from \cite[Theorem 2.3.2 p.112]{mazya} that
$$
\NRME{f_+^2}{\Phi} \leq
4 B(\Phi, \{f \leq c\}) \int |\nabla f_+|^2 d \nu ,
$$
where $B(\Phi,\{f \leq c\} )$ is the smallest constant so that for every
$A \subset \{f \leq c\}$,
$$
\NRME{\ind_A}{\Phi} \leq B(\Phi,\{f \leq c\} )
\mathrm{Cap}_\nu(A,\{f \leq c\}) .
$$
A similar result holds for $f_-$. Thus, by definition of
$\mathrm{Cap}_\nu(A,\mu)$ and $B(\Phi)$,
we get $B(\Phi,\{f \leq c\} ) \leq B(\Phi)$ and in turn
\begin{eqnarray*}
\NRME{(f-\mu(f))^2}{\Phi} & \leq &
8(1+k) B(\Phi)
\left( \int |\nabla f_+|^2 d \nu + \int |\nabla f_-|^2 d \nu \right) \\
& \leq &
8(1+k) B(\Phi)   \int |\nabla f|^2 d \nu .
\end{eqnarray*}
In the last inequality we used that, since $f$ is locally Lipschitz and
$\nu$ is absolutely continuous, the set $\{f=c\} \cap \{\nabla f \neq 0 \}$
is $\nu$-negligible. This proves the first part of the criterium.

For the other part, take a Borel set $A \subset \dR^n$ with
$\mu(A) < \frac 12$ and a function $f$ such that
$\mu(\{f=0\}) \geq \frac 12$ and $ \ind_{\{f \neq 0\}} \geq f \geq \ind_A$.
Set ${\cal G}=\{g: \dR^n \rightarrow \dR ; \int \Phi^*(g) d\mu \leq 1 \}$
where $\Phi^*$ is the conjugate function of $\Phi$.
By \eqref{in:eqn} we have
\begin{eqnarray*}
2 \NRME{(f-\mu(f))^2}{\Phi}
&\geq&
\sup_{ g \in \cal G} \int((f-\mu(f))^2 |g| d\mu
\geq
\sup_{ g \in \cal G} \int_A((f-\mu(f))^2 |g| d\mu \\
& = &
 (1-\mu(f))^2 \sup_{ g \in \cal G} \int_A|g| d\mu
\geq
 (1-\mu(f))^2 \NRME{\ind_A}{\Phi} .
\end{eqnarray*}
Since $f \leq \ind$, we get
$\mu(f) \leq \mu(f \neq 0) \leq \frac 12$.
Thus $(1-\mu(f))^2  \geq \frac 1{4}$.
It follows that
$$
\frac1{8} \NRME{\ind_A}{\Phi} \leq \NRME{(f-\mu(f))^2}{\Phi}
\leq C_\Phi \int |\nabla f|^2 d \nu .
$$
The result follows by definition of the capacity. This ends the proof
in any finite dimension.

Consider now $n=1$. Let $m$ be a median of $\mu$ and define
$f_+=(f-f(m)) \ind_{(m,+\infty)}$ and
$f_-=(f-f(m)) \ind_{(-\infty,m)}$. Note that $(f_+ + f_-)^2=f_+^2 + f_-^2$.
By our assumption and a similar computation as in the general case,
\begin{eqnarray*}
\NRME{(f-\mu(f))^2}{\Phi}
& = &
\NRME{(f-f(m) - \mu(f-f(m)))^2}{\Phi} \\
&\leq&
2 (1+k) \NRME{(f-f(m))^2}{\Phi}
=
2(1+k) \NRME{(f_+ + f_-)^2}{\Phi} \\
& \leq &
2(1+k)(\NRME{f_+^2}{\Phi}+\NRME{f_-^2}{\Phi}) .
\end{eqnarray*}
From \cite[Proposition 2]{barthe-roberto} (which originally comes from
\cite{bobkov-gotze}, see also \cite{chen}), it follows that
$$
\NRME{f_+^2}{\Phi} \leq 4 B_+(\Phi) \int_m^\infty {f_+'}^2 d\nu .
$$
Since a similar bound holds for $f_-$, summing up we get that
$C_\Phi \leq 8(1+k) \max(B_+(\Phi),B_-(\Phi))$.

Next, fix $x > m$ and consider the following function defined on the real line
$$
h(y)=\left\{
\begin{array}{ll}
0 & \mbox{for } y \leq m \\
\int_m^y \frac{1}{\rho_\nu} & \mbox{for } m \leq y \leq x \\
\int_m^x \frac{1}{\rho_\nu} & \mbox{for } y \geq x .
\end{array}
\right.
$$
Starting as previously, we get that
\begin{eqnarray*}
2 \NRME{(h-\mu(h))^2}{\Phi}
& \geq &
\sup_{ g \in \cal G} \int_{[x,\infty)} ((h-\mu(h))^2 |g| d\mu \\
& \geq &
\left( \int_m^x \frac{1}{\rho_\nu} - \mu(h)\right)^2 \NRME{\ind_{[x,\infty)}}{\Phi} .
\end{eqnarray*}
Then, since $x>m$ and $h \leq \int_m^x \frac{1}{\rho_\nu}$,
$$
\mu(h) \leq  \mu((m,\infty)) \int_m^x \frac{1}{\rho_\nu}
\leq \frac 12 \int_m^x \frac{1}{\rho_\nu} .
$$
Therefore, $\int_m^x \frac{1}{\rho_\nu} - \mu(h)
\geq \frac 12 \int_m^x \frac{1}{\rho_\nu}$.
Applying the Orlicz-Sobolev inequality to this special function $h$,
we get
$$
\frac 14 \left( \int_m^x \frac{1}{\rho_\nu} \right)^2
\NRME{\ind_{[x,\infty)}}{\Phi}
\leq 2 \NRME{(h-\mu(h))^2}{\Phi}
\leq 2 C_\Phi \int {h'}^2 d\nu
= 2 C_\Phi \int_m^x \frac{1}{\rho_\nu} .
$$
This gives for any $x >m$,
$$
\NRME{\ind_{[x,\infty)}}{\Phi}
\int_m^x \frac{1}{\rho_\nu}  \leq 8 C_\Phi .
$$
The same bound holds for $x < m$ and the result follows by definition
of $B_+(\Phi)$ and $B_-(\Phi)$.
\end{proof}

The explicit criterium in dimension 1 leads to the following result.

\begin{proposition} \label{prop:c}
Let $\Phi$ be an Young function and fix $k \in (0,+\infty)$
such that
$\NRME{\mu(f)^2}{\Phi} \leq k \NRME{f^2}{\Phi}$, for any function $f$ with $f^2 \in \L_\Phi (\mu)$.
Let $V : \dR \rightarrow \dR$ be a $\C{1}$ function such that
$d\mu(x)=e^{-V(x)}dx$ is a probability measure.
Furthermore assume that\\
\noindent $(i)$ there exists a constant $A>0$ such that for
$|x| \geq A$, $V$ is $\C{2}$ and $\mathrm{sign}(x) V'(x) > 0$,\\
\noindent $(ii)$ $\lim_{|x| \rightarrow  \infty} \frac{V''(x)}{V'(x)^2} = 0$,\\
\noindent $(iii)$
$\liminf_{|x| \rightarrow \infty} V'(x)e^{-V(x)} \Phi^{-1}(V'(x) e^{V(x)}) >0$.
\\
Then there exists
a constant $C_\Phi$ (that may depend on $k$)
such that for every smooth function $f : \dR \rightarrow \dR$,
one has $$
\NRME{(f-\mu(f))^2}{\Phi} \leq C_\Phi \int {f'}^2 d\mu .
$$
\end{proposition}

\begin{proof}
The proof is similar to \cite[Theorem 6.4.3 (Chapter 6)]{ane}. Let $m$ be
a median of $\mu$. Under assumptions $(i)$ and $(ii)$, when $x$ tends to
infinity, one has (see {\it e.g.} \cite[Chapter 6]{ane})
\begin{eqnarray*}
\int_m^x e^{V(t)}dt \sim \frac{e^{V(x)}}{V'(x)}
\qquad
\mbox{and}
\qquad
\int_x^\infty e^{-V(t)}dt \sim \frac{e^{-V(x)}}{V'(x)} .
\end{eqnarray*}
Thus, for $x>m$,
\begin{eqnarray*}
\NRME{\ind_{[x,\infty)}}{\Phi} \int_m^x e^{V(t)}dt
& = &
\frac{1}{\Phi^{-1}(1/\mu([x,\infty)))} \int_m^x e^{V(t)}dt \\
& \sim &
\frac{1}{ V'(x)e^{-V(x)} \Phi^{-1}(V'(x) e^{V(x)})} .
\end{eqnarray*}
By hypothesis $(iii)$ this quantity is bounded on $[A',\infty)$ for some
$A' \geq m$. Since it is continuous on $[m,A']$, it is bounded on $(m,\infty)$.
It follows that $B_+(\Phi)$ and $B_-(\Phi)$, (defined in Theorem
\ref{th:cos}), are bounded. We conclude with Theorem \ref{th:cos}.
\end{proof}

In general the capacity can be difficult to compute. However it provides a nice
interfacing tool to prove equivalences between inequalities.
Indeed, a criterium involving capacity also holds for
general Beckner-type inequalities as we will see now.
The two general criterium will allows us to prove an equivalence between the Orlicz-Sobolev
inequality and the Beckner-type inequalities.
Our main motivation here is that the latter naturally tensorises. Thus,
dimension free Orlicz Sobolev inequalities will follow from Beckner-type inequality.

Combining Theorem 9 and Lemma 8
of \cite{barthe-cattiaux-roberto} we get the following:

\begin{theorem}[\cite{barthe-cattiaux-roberto}]\label{th:cgb}
Let $T:[0,1] \rightarrow \dR^+$ be
non-decreasing and such that $x \mapsto T(x)/x$
is non-increasing. Let $\mu $ and $\nu$ be
two absolutely continuous measures on $\dR^n$ with $\mu(\dR^n)=1$.
Let $C_T$ be the optimal constant such that for
every smooth $f: \dR^n \rightarrow \dR$ one has
\begin{equation}\label{in:gb}
\sup_{p \in (1,2)} \frac{\int f^2 d\mu - (\int |f|^p d\mu)^\frac2p}{T(2-p)}
\leq C_T \int |\nabla f|^2 d \nu .
\end{equation}
Then, $\frac 16 B(T) \leq C_T \leq 20 B(T)$, where $B(T)$
is the smallest constant so that
every Borel set $A \subset \dR^n$ with $\mu(A) < \frac 12$ satisfies
$$
\frac{\mu(A)}{T\bigl(1/\log(1+\frac{1}{\mu(A)})\bigr)}
\leq B(T) \mathrm{Cap}_\nu(A,\mu) .
$$
\end {theorem}

Now, using the previous two theorems,
one can see that the Orlicz-Sobolev inequality
({\textbf{O-S}})
is equivalent, up to universal constant,
to the general Beckner-type inequality
\eqref{in:gb}.

\begin{corollary}\label{cor:e}
Let $\mu$ and $\nu$ be two absolutely continuous measures on $\dR^n$
with $\mu(\dR^n)=1$.
Let $T:[0,1] \rightarrow \dR^+$ be non-decreasing and such that $x \mapsto T(x)/x$
is non-increasing. Denote by $C_T$ the optimal constant such that
for every smooth $f: \dR^n \rightarrow \dR$ one has
\begin{equation*}
\sup_{p \in (1,2)} \frac{\int f^2 d\mu - (\int |f|^p d\mu)^\frac2p}{T(2-p)}
\leq C_T  \int |\nabla f|^2 d \nu .
\end{equation*}
Let $\Phi$ be a Young function and let $k \in (0, + \infty)$
be such that for any function $f$ with $f^2 \in \L_\Phi (\mu)$,
$\NRME{\mu(f)^2}{\Phi} \leq k \NRME{f^2}{\Phi}$. Let
$C_\Phi$ the optimal constant such that for every smooth
$f : \dR^n \rightarrow \dR$ one has
$$
\NRME{(f-\mu(f))^2}{\Phi} \leq C_\Phi \int |\nabla f|^2 d \nu .
$$
Finally, assume that there exists two positive constant $c_1$ and $c_2$ such that
$$
c_1 x T\bigl(\frac{1}{\log(1+x)}\bigr)
\leq \Phi^{-1}(x) \leq
c_2  x T\bigl(\frac{1}{\log(1+x)}\bigr)  \qquad \forall x > 2 .
$$
Then,
$$
\frac{c_1}{48(1+k)}\ C_\Phi \leq C_T \leq 160 c_2\ C_\Phi .
$$
\end{corollary}

\begin{proof}
The last assumption on $T$ and $\Phi^{-1}$ is equivalent to
$$
\frac{1}{c_2} \frac{y}{T\bigl(1/\log(1+\frac{1}{y})\bigr)} \leq
\frac{1}{\Phi^{-1}(1/y)} \leq
\frac{1}{c_1} \frac{y}{T\bigl(1/\log(1+\frac{1}{y})\bigr)} \qquad
\forall y \in (0, \frac 12) .
$$
Since $\NRME{\ind_A}{\Phi}= \frac{1}{\Phi^{-1}(1/\mu(A))}$, it follows that
$\frac{1}{c_2} B(T) \leq B(\Phi) \leq \frac{1}{c_1} B(T)$, where
$B(\Phi)$ and  $B(T)$ are defined in Theorem \ref{th:cos} and \ref{th:cgb} respectively.
The result follows from Theorem \ref{th:cos} and \ref{th:cgb}.
\end{proof}

\subsection*{Example: $T_\beta(x)=|x|^\beta$} 
An important example is given by $T_\beta(x)=|x|^\beta$ with $\beta \in [0,1]$.
This correspond to the Lata{\l}a and Oleszkiewicz inequality (in short L-O inequality)
\cite{latala-Oleszkiewicz}.

Consider the Young function $\Phi_\beta (x)= |x| [\log(1+|x|)]^\beta$ with
$\beta \in [0,1]$. Then, we claim that
\begin{equation} \label{eq:inverse}
\frac{y}{[\log(1+y)]^\beta}
\leq
\Phi_\beta^{-1}(y)
\leq
2 \frac{y}{[\log(1+y)]^\beta}
\qquad \forall y >2 .
\end{equation}
Indeed,
$$
\Phi_\beta( \frac{x}{[\log(1+x)]^\beta} ) =
x \frac{[\log(1+x (\log(1+x))^{-\beta})]^\beta}{[\log(1+x)]^\beta}
\qquad \forall x \geq 0 .
$$
Note that for $x \geq e -1$, $1+x(\log(1+x))^{-\beta} \leq 1+x$.
This leads to
$\Phi_\beta( \frac{x}{[\log(1+x)]^\beta} ) \leq x$
for $x >2$. The first inequality in \eqref{eq:inverse} follows.

On the other hand, it is not difficult to check that for any
$\gamma \in [0,1]$, any $x \geq e-1$,
$$
1+x(\log(1+x))^{-\beta}
\geq
\frac{1+x}{[\log(1+x)]^\beta}
\geq \left( \frac{e(1-\gamma)}{\beta} \right)^\beta (1+x)^\gamma.
$$
It follows for $\gamma=1-(\beta/e)$ that
$\Phi_\beta( \frac{x}{[\log(1+x)]^\beta} ) \geq \gamma^\beta x
\geq \frac{e-1}{e}x$.
Thus, for any $y \geq (e-1)^2/e \simeq 1.09$,
$$
\Phi_\beta^{-1}(y)
\leq
\frac{e}{e-1}\frac{y}{\log(1+\frac{ey}{e-1})^\beta}
\leq
\frac{e}{e-1}\frac{y}{[\log(1+y)]^\beta} .
$$
The result follows.

We are in position to prove a family of Orlicz-Sobolev inequalities.

\begin{corollary}\label{cor:alpha}
Let $\alpha \in [1,2]$,
$\beta = 2(1-\frac1\alpha) \in [0,1]$
and $\Phi_\beta (x)= |x| [\log(1+|x|)]^\beta$.
Then, for any integer $n$, the probability measure
on $\dR^n$,
$d\mu_\alpha^n(x) =
Z_\alpha^{-n} \exp \{-\sum_{i=1}^n |x_i|^\alpha \} dx$
satisfies for any smooth function $f: \dR^n \rightarrow \dR$,
\begin{equation}\label{eq:oslo}
\NRME{(f-\mu_\alpha^n(f))^2}{\Phi_\beta}
\leq
C \int |\nabla f|^2 d\mu_\alpha^n
\end{equation}
for some universal constant $C$ independent of $n$ and $\alpha$.
\end{corollary}

\begin{proof}
Fix an integer $n$, $\alpha \in [1,2]$,
the corresponding $\beta$ and let $T_\beta(x)=|x|^\beta$.
It is proved in \cite{latala-Oleszkiewicz} that
the measure $d\mu_\alpha^n$ on $\dR^n$ satisfies the general
Beckner type inequality \eqref{in:gb} with $T=T_\beta$
and constant $C(T_\beta)$ independent of $n$ and $\alpha$
(for the uniformity in $\alpha$, see
\cite[section 7]{barthe-cattiaux-roberto}).
Then the result follows by our previous claim (inequality \eqref{eq:inverse})
on  $\Phi_\beta^{-1}$ and Corollary \ref{cor:e}
(it is easy to check that
$\NRME{\mu_\alpha(f)^2}{\Phi_\beta} \leq e \NRME{f^2}{\Phi_\beta}$
from Remark \ref{remark_norm_mu_f2}).
\end{proof}

\begin{remark}
The family of inequalities in Corollary \ref{cor:alpha} is an
interpolation family between Poincar{\'e}, for $\Phi(x)=|x|$, and
the logarithmic Sobolev inequality, for $\Phi(x)=|x| \log(1+|x|)$
(see \cite{bobkov-gotze}).
\end{remark}

\begin{remark}

To prove that inequality \eqref{eq:oslo} holds in dimension 1,
we could have used
Proposition \ref{prop:c} together with \eqref{eq:inverse}. Moreover,
given $\beta \in [0,1]$,
Proposition \ref{prop:c} insures that \eqref{eq:oslo} holds
for any $\alpha \geq \alpha(\beta)$ where $\beta = 2(1-\frac1{\alpha(\beta)})$
and does not hold for  $\alpha < \alpha(\beta)$.
\end{remark}


\subsection*{L-O inequality for Gibbs measures}

The following result provides a precise asymptotic of the coefficient
in L-O inequality as well as plays a vital role in a construction of
examples non-product measures satisfying this inequality.

\begin{theorem}
$(i)$ Let $p\in[1,2]$. Then,
$$
\NRM{f}_2^2 - \NRM{f}_p^2
\leq
(p-1) \left( \NRM{f-\mu (f)}_2^2 - \NRM{f-\mu (f)}_p^2 \right)
+ (2-p) \NRM{f-\mu (f)}_2^2 .
$$
Hence, if with some $C\in(0,\infty)$ and $\beta\in(0,1)$
$$
\NRM{f-\mu (f)}_2^2 - \NRM{f-\mu (f)}_p^2
\leq
C (2-p)^{\beta} \NRM{\nabla f}_2^2 ,
$$
and for some $M\in(0,\infty)$
$$
M \NRM{f-\mu (f)}_2^2 \leq  \NRM{\nabla f}_2^2 ,
$$
then
$$
\NRM{f}_2^2 - \NRM{f}_p^2
\leq
\left((p-1) C (2-p)^{\beta}+ (2-p)M\right)\NRM{\nabla f}_2^2 .
$$
$(ii)$ \emph{(Mild Perturbation Lemma)}
Suppose $\nu$ satisfies the following L-O inequality
$$
\NRM{f}_{\L_2(\nu)}^2 - \NRM{f}_{\L_p(\nu)}^2
\leq
C (2-p)^{\beta}
\NRM{\nabla f}_{\L_2(\nu)}^2
$$
and let $d\mu = \rho d\nu$ with
$\delta U\equiv \sup (\log \rho) - \inf(\log\rho) < \infty$.
Then
$$
\NRM{f}_{\L_2(\mu)}^2 - \NRM{f}_{\L_p(\mu)}^2
\leq
e^{\delta U} C (2-p)^{\beta} \NRM{\nabla f}_{\L_2(\mu)}^2 .
$$
\end{theorem}

\begin{proof} (\cite{wang})
$(i)$ The first inequality follows from the following convexity
property of the $\L_p(\mu)$ norm for $p\in[1,2]$
$$
\NRM{f}_{\L_p(\mu)}^2
\geq
\mu (f)^2 +(p-1)\NRM{f-\mu (f)}_{\L_p(\mu)}^2 ,
$$
see e.g. \cite{bcl}.
(sese also \cite{wang}, \cite[Lemma 8]{barthe-roberto}).
This together with spectral gap inequality and L-O inequality for
$f-\mu (f)$ imply the L-O inequality for $f$ with the improved coefficient.\\
$(ii)$ We note first that for $p\in(1,2)$, with $A\equiv \frac {2-p}2
\big(\frac p2\big)^{\frac p{2-p}}$, we have
$$
\NRM{f}_{\L_2(\mu)}^2 -\NRM{f}_{\L_p(\mu)}^2 =
\inf_{t>0} \mu\left( f^2 - t |f|^p + At^{\frac 2{2-p}} \right) .
$$
Since by Young inequality
$$
z^pt = \Big[\Big(\frac2p\Big)^{\frac p2}z^{p}\Big]
\cdot\Big[\Big(\frac2p\Big)^{-\frac p2}t\Big]
\leq z^2 + \frac{2-p}2 \Big(\frac2p\Big)^{-\frac p{2-p}}t^{\frac 2{2-p}}
= z^2 + A t^{\frac 2{2-p}}
$$
the integrand in the above is nonnegative. Hence, if $d\mu =\rho d\nu$,
we get
\begin{eqnarray*}
\inf_{t>0} \mu\left( f^2 - t |f|^p + At^{\frac p{2-p}} \right)
&\leq &
\sup(\rho)\inf_{t>0} \nu\left( f^2 - t |f|^p + At^{\frac p{2-p}} \right) \\
& \leq &
\sup(\rho) C(2-p)^\beta \nu \left(|\nabla f|^2 \right) \\
& \leq &
\frac{\sup(\rho)}{\inf(\rho)} C(2-p)^\beta \mu \left(|\nabla f|^2 \right).
\end{eqnarray*}
\end{proof}

Starting from the product measure satisfying L-O inequality,
using the Mild Perturbation Lemma we see that one can construct
a local specification
for which each finite volume conditional expectation $E_\Lambda$
(defined as a mild perturbation of the product measure), satisfies this inequality.
This together with the suitable conditioning expansion
based on the following step
\begin{eqnarray*}
{\boldsymbol\mu}(f^2) - \left({\boldsymbol\mu}(f^p)\right)^\frac2p
& = &
{\boldsymbol\mu}\left(E_\Lambda (f^2) -
\left(E_\Lambda (f^p)\right)^\frac2p\right) \\
&& \quad +
{\boldsymbol\mu} \left(\left[E_\Lambda (f^p)^\frac1p\right]^2\right)
- {\boldsymbol\mu}\left(\left[E_\Lambda(f^p)^\frac1p\right]^p\right)^{\frac2p}
\end{eqnarray*}
under suitable mixing condition (the same as the one used in the case of log-Sobolev inequality), allows to prove the following result
(see \cite{GZ03} for details).

\begin{theorem}
Suppose a local specification is mixing and satisfies {\textbf{L-O}}
inequality with the index $\beta\in(0,1)$.
Then the corresponding Gibbs measure ${\boldsymbol\mu}$ satisfies
$$
{\boldsymbol\mu}\left(f^2\right) - {\boldsymbol\mu}\left(f^p\right)^\frac2p
\leq C(2-p)^\beta {\boldsymbol\mu}\left(|\nabla f|^2\right)
$$
with a constant $C\in(0,\infty)$ independent of a function $f$.
\end{theorem}


\section{O-S inequality and Decay to equilibrium}

In this section we prove that the semi-group naturally
associated to a measure $\mu$
satisfying an Orlicz-Sobolev inequality decays exponentially
fast in $\L_\Phi(\mu)$. This
result is new and strengthens a well know fact for the Poincar{\'e}
inequality and decay in $\L_2(\mu)$.
We start with a modified Orlicz-Sobolev inequality.
As before, throughout below we consider the following setup.
Let $d\mu(x)=e^{V(x)}dx$ be a probability
measure on $\dR^n$ associated to the differentiable potential $V$.
Let $\GI=\Delta - \nabla V \cdot \nabla$ be a symmetric in $\L_2(\mu)$
diffusion generator
and $(\PT{t})_{t \geq 0}$ its associated semi-group.

\begin{theorem}\label{th:decay}
Consider a Young function $\Phi$
satisfying $x \Phi'(x) \leq B \Phi(x)$ for every $x$ and
some constant $B$.
Then, the following are equivalent\\
\noindent $(i)$ There exists a constant $C_\Phi$ such that for any
smooth function $f: \dR^n \rightarrow \dR$,
\begin{equation} \label{in:mos}
\NRM{f-\mu(f)}_\Phi^2 \leq
C_\Phi \int |\nabla f|^2 \Phi'' \left( \frac{f-\mu(f)}{\NRME{f-\mu(f)}{\Phi}} \right) d \mu .
\end{equation}
\noindent $(ii)$  There exists a constant $M\in(0,\infty)$  such that for any
smooth function $f: \dR^n \rightarrow \dR$, for
any $t \geq 0$,
$$
\NRM{\PT{t}f-\mu(f)}_\Phi^2 \leq e^{-Mt} \NRM{f-\mu(f)}_\Phi^2 .
$$
Furthermore, $(i)$ implies $(ii)$ with $M=2/(BC_\Phi)$, and  $(ii)$ implies $(i)$
with $C_\Phi=2/M$.
\end{theorem}

\begin{remark}
Note that if $\Phi$ satisfies the $\Delta_2$-condition $\Phi(2x) \leq C \Phi(x)$ for every $x$,
then
$$
x \Phi'(x) \leq \int_x^{2x} \Phi'(t) dt = \Phi(2x) -\Phi(x) \leq (C-1) \Phi(x) .
$$
Thus the condition on the Young function $\Phi$
is satisfied as soon as the $\Delta_2$-condition
is satisfied.
\end{remark}

\begin{proof}
Without loss of generality for a smooth non zero function $f$, we can assume that $\mu(f) = 0$.
Let $N(t) = \NRME{\PT{t}f}{\Phi}$.

By definition of the Luxembourg norm, we have
$\int\Phi \left( \frac{\PT{t}f}{N(t)}  \right) d\mu = 1$.
A differentiation and the chain rule
formula $\int \Phi'(g) L g d\mu = - \int \Phi''(g) |\nabla g|^2 d\mu$ give
\begin{eqnarray*}
\frac{N'(t)}{N(t)} \int \Phi' \left( \frac{\PT{t}f}{N(t)} \right) \frac{\PT{t}f}{N(t)} d\mu
& = &
\int \frac{\GI \PT{t}f}{N(t)} \Phi' \left( \frac{\PT{t}f}{N(t)} \right)  d\mu \\
& = &
-\frac{1}{N^2(t)}  \int \Phi'' \left( \frac{\PT{t}f}{N(t)} \right) |\nabla \PT{t}f|^2 d\mu .
\end{eqnarray*}
We will first show that $(i) \Rightarrow (ii)$.
Since $\Phi$ is a Young function,
it is convex and for any $x$, $x \Phi'(x) \geq 0$.
It follows at first that $N'(t) \leq 0$. Furthermore, by hypothesis
$x \Phi'(x) \leq B \Phi(x)$. Thus, using the property that
$\int\Phi \left( \frac{\PT{t}f}{N(t)}  \right) d\mu = 1$, we get by $(i)$ that
$$
B \frac{N'(t)}{N(t)} \leq -\frac{1}{N^2(t)}
\int \Phi'' \left( \frac{\PT{t}f}{N(t)} \right) |\nabla \PT{t}f|^2 d\mu
\leq
- \frac{1}{N^2(t)}
\frac{1}{C_\Phi} N^2(t) = - \frac{1}{C_\Phi},
$$
which gives the expected result.

Now we show that $(ii) \Rightarrow (i)$. Let
$u(t)=e^{Mt}\NRM{\PT{t}f-\mu(f)}_\Phi^2$.
Point $(ii)$ exactly means that $u'(t) \leq 0$. Hence
$ Me^{Mt} N^2(t) + 2e^{Mt} N'(t) N(t) \leq 0$ which leads to
\begin{eqnarray*}
M\ N^2(t) \leq - 2 N'(t) N(t)
& = &
2 \frac{\int \Phi'' \left( \frac{\PT{t}f}{N(t)} \right)
|\nabla \PT{t}f|^2 d\mu }
{ \int \Phi' \left( \frac{\PT{t}f}{N(t)} \right) \frac{\PT{t}f}{N(t)} d\mu} \\
& \leq &
2 \int \Phi'' \left( \frac{\PT{t}f}{N(t)} \right) |\nabla \PT{t}f|^2 d\mu .
\end{eqnarray*}
In the last inequality we used the fact that,
since $\Phi$ is convex and $\Phi(0)=0$,
for every $x$, $x \Phi'(x) \geq \Phi(x)$ and
$\int\Phi \left( \frac{\PT{t}f}{N(t)}  \right) d\mu = 1$.
The latter inequality applied at $t=0$ gives
the expected result. This ends the proof.
\end{proof}

\begin{remark}
When $\Phi(x)=x^2$, $\NRM{f}_\Phi^2 = \NRM{f}_2^2$ and $\Phi''(x)=2$.
Thus Theorem \ref{th:decay} recover the well known equivalence between the exponential
decay of the semi-group in $\L_2$-norm and the Poincar{\'e} inequality.
\end{remark}

The behavior of $\Phi''$ seems to play an important role.
In particular, under additional strict positivity assumption
we prove the following result involving the Orlicz-Sobolev inequalities.

\begin{corollary}\label{cor:ex}
Consider a Young function $\Phi$ and set $\Phi_2(x)=\Phi(x^2)$. Assume that
 $x \Phi_2'(x) \leq B \Phi_2(x)$ for every $x$ and
some constant $B$, and $\Phi_2'' \geq \ell >0$.
Assume that there exists a constant $C_\Phi$ such that
for any smooth function $f: \dR^n \rightarrow \dR$,
$$
\NRM{(f-\mu(f))^2}_\Phi \leq C_\Phi \int |\nabla f|^2  d \mu .
$$
Then, for any smooth function $f$, for any $t \geq 0$,
$$
\NRM{(\PT{t}f-\mu(f))^2}_\Phi \leq e^{-Mt}
\NRM{(f-\mu(f))^2}_\Phi .
$$
with $M = \frac{2 \ell}{BC_\Phi}$
\end{corollary}

\begin{proof}
It is enough to check that for any function $f$, we have
\begin{eqnarray*}
\NRM{f-\mu(f)}_{\Phi_2}^2 & = & \NRM{(f-\mu(f))^2}_\Phi  \leq
C_\Phi \int |\nabla f|^2  d \mu  \\
& \leq & \frac{C_\Phi}{\ell}
\int |\nabla f|^2
\Phi_2'' \left( \frac{f-\mu(f)}{\NRM{f-\mu(f)}_{\Phi}} \right) d \mu ,
\end{eqnarray*}
and to apply Theorem \ref{th:decay}.
\end{proof}

In Corollary \ref{cor:alpha} we proved that a family of Orlicz-Sobolev
inequalities hold for $\Phi_\beta (x)= |x| [\log(1+|x|)]^\beta$,
$\beta \in [0,1]$.
Actually we cannot apply the previous result to this family of norms,
simply because
$\Phi_{\beta,2}''(0)=0$ and thus there is no bound of
the type  $\Phi_{\beta,2}'' \geq \ell >0$
(here  $\Phi_{\beta,2} (x)= x^2 \log(1+x^2)^\beta$).
However we can get rid of this problem by means
of equivalence of norms. 

\begin{proposition} \label{prop:3c}
Let $\alpha \in [1,2]$, $\beta = 2(1-\frac1\alpha) \in [0,1]$
and for $\gamma \geq 1$,
$\Phi_\beta^\gamma (x)= |x| \log(\gamma+|x|)^\beta$.
Let $d\mu_\alpha^n(x) = Z_\alpha^{-n} \exp \{-\sum_{i=1}^n |x_i|^\alpha \}dx$
be a probability measure on $\dR^n$,
$\GI=\Delta + \nabla V \cdot \nabla$ with
$V=\sum_{i=1}^n |x_i|^\alpha$
be a symmetric (in $\L_2(\mu_\alpha)$) diffusion generator
and $(\PT{t})_{t \geq 0}$ its associated semi-group.
Let $C$ be the coefficient appearing in the Orlicz-Sobolev inequality of Corollary \ref{cor:alpha}.

Then, for any $\gamma > 1$, any $\beta \in [0,1]$, any integer $n$,
any function $f$ and any $t \geq 0$,
$$
\NRM{(\PT{t}f-\mu_\alpha^n(f))^2}_{\Phi_\beta^\gamma} \leq
e^{-c_1t}
\NRM{(f-\mu_\alpha^n(f))^2}_{\Phi_\beta^\gamma} ,
$$
with
$c_1=\frac{(\log \gamma)^\beta}{4^\beta C (1+e (\log \gamma)^\beta)}$.\\
While for any $\beta \in [0,1]$, any integer $n$,
any function $f$ and any $t \geq 0$,
$$
\NRM{(\PT{t}f-\mu_\alpha^n(f))^2}_{\Phi_\beta^1}
\leq
\left\{ \begin{array}{ll}
\NRM{(f-\mu_\alpha^n(f))^2}_{\Phi_\beta^1} &
\mbox{ for } t \leq 4^\beta C e \\
\frac{t}{4^\beta C}e^{-\frac{t}{4^\beta C e}}
\NRM{(f-\mu_\alpha^n(f))^2}_{\Phi_\beta^1} & \mbox{ for } t \geq 4^\beta C e .
\end{array} \right.
$$
\end{proposition}

\begin{proof}
Fix $\gamma > 1$, an integer $n$ and $\beta \in [0,1]$.
Then note that from the equivalence of Orlicz norms corresponding
for different $\gamma$, (see point $(i)$ of Lemma \ref{lem:eq}
below (with $\gamma=1$)),
and Corollary \ref{cor:alpha}, for any sufficiently smooth function $f$, one has
\begin{eqnarray*}
\NRME{(f-\mu_\alpha^n(f))^2}{\Phi_\beta^\gamma}
& \leq &
\left(1+e(\log \gamma)^\beta \right)
\NRME{(f-\mu_\alpha^n(f))^2}{\Phi_\beta^1}\\
& \leq &
C \left(1+e (\log \gamma)^\beta \right) \int |\nabla f|^2 d \mu_\alpha^n .
\end{eqnarray*}
On the other hand, by point $(ii)$ and $(iii)$ of Lemma \ref{lem:eq} below
we can apply Corollary \ref{cor:ex} with $B=4^{1+\beta}$ and
$\ell=2(\log \gamma)^\beta$. It follows that
$$
\NRM{(\PT{t}f-\mu_\alpha^n(f))^2}_{\Phi_\beta^\gamma} \leq
\exp\left\{-
\frac{(\log \gamma)^\beta}{4^\beta C (1+e (\log \gamma)^\beta)} t
\right\}
\NRM{(f-\mu_\alpha^n(f))^2}_{\Phi_\beta^\gamma}
$$
which gives the first part of the result.

For the second part, we use twice the latter inequality together
with point $(i)$ of Lemma \ref{lem:eq} below to get for any
$\gamma \geq 1$,
\begin{eqnarray*}
\NRM{(\PT{t}f-\mu_\alpha^n(f))^2}_{\Phi_\beta^1}
& \leq &
\NRM{(\PT{t}f-\mu_\alpha^n(f))^2}_{\Phi_\beta^\gamma} \leq  e^{-c_1t}
\NRM{(f-\mu_\alpha^n(f))^2}_{\Phi_\beta^\gamma} \\
& \leq &
\left(1+e (\log \gamma)^\beta \right)  e^{-c_1t}
\NRM{(f-\mu_\alpha^n(f))^2}_{\Phi_\beta^1}
\end{eqnarray*}
with $c_1=\frac{(\log \gamma)^\beta}{4^\beta C (1+e (\log \gamma)^\beta)}$.
The result follows from an optimization over $\gamma \geq 1$
and the decrease of $N(t)$ proved before.
\end{proof}

\begin{lemma}\label{lem:eq} 
For $\beta \in [0,1]$ and $\gamma \geq 1$,
let $\Phi_\beta^\gamma (x)= |x| \log(\gamma+|x|)^\beta$ and
$\Phi_{\beta,2}^\gamma (x) = \Phi_\beta^\gamma (x^2)$. Then,\\
\noindent $(i)$ for any $1 \leq \gamma \leq \gamma'$,
$$
\NRM{\cdot}_{\Phi_{\beta,2}^\gamma} \leq
\NRM{\cdot}_{\Phi_{\beta,2}^{\gamma'}} \leq C_{\gamma ,\gamma'}
\NRM{\cdot}_{\Phi_{\beta,2}^\gamma} .
$$
with $C_{\gamma ,\gamma'} \equiv
\bigl[1+(1+(e-\gamma)_+)
\bigl(\log \frac{\gamma'}{\gamma}\bigr)^\beta\bigr]^\frac{1}{2}$, where  $(x)_+:=\max(x,0)$.\\
\noindent $(ii)$ For any $x$,
${\Phi_{\beta,2}^\gamma} ''(x) \geq 2(\log \gamma)^\beta$.\\
\noindent $(iii)$ For any $x$,
$\Phi_{\beta,2}^\gamma (2x) \leq 4^{1+\beta} \Phi_{\beta,2}^\gamma (x)$.
\end{lemma}

\begin{proof}
First, $(i)$ follows from Lemma \ref{lemma_injection_L1},
provided in the Appendix, since for any $1 \leq \gamma \leq \gamma'$ one has
\begin{align*}
\Phi_{\beta}^{\gamma'} (x) = |x| \, \left(\log \frac{\gamma'}{\gamma}
 + \log \left(\gamma + \frac{\gamma}{\gamma'} |x| \right)
\right)^\beta
 \leq \left(\log \frac{\gamma'}{\gamma}\right)^\beta \, |x|
+ \Phi_{\beta}^{\gamma} (x).
\end{align*}
We also made use of the bound \eqref{norm_injection}
for $\Phi =\Phi_{\beta}^{\gamma}$,  $\tau =1$ and $M = (e-\gamma)_+$.

Now, we may easily check that ${\Phi_{\beta,2}^\gamma} ''$ is
non-decreasing and thus
greater than ${\Phi_{\beta,2}^\gamma} ''(0)=2(\log \gamma)^\beta$.
This gives $(ii)$.

Using $\gamma + 4 x^2 \leq (\gamma + x^2)^4$ (recall that
$\gamma \geq 1$), we get
$$
\Phi_{\beta,2}^\gamma (2x)
 =
4 x^2 (\log (\gamma + 4 x^2))^\beta
 \leq
4 x^2 (\log (\gamma +  x^2)^4)^\beta
 =
4^{1+\beta} \Phi_{\beta,2}^\gamma (x) .
$$
The proof is complete.
\end{proof}

\subsection*{Monotone Functionals}  \label{Monotone Functionals}

In Proposition \ref{prop:3c} the semi-group is not decaying exponentially
to equilibrium in particular for $\Phi_1^1$. We shall see in this section
that a modification (a time-averaging) of the functional will
satisfies an  exponential decay.

The following inequality was shown in \cite[Proposition 4.1]{bobkov-gotze}
$$
\frac{2}{3} \NRM{(f-\mu_2(f))^2}_{\Phi_1^1}
\leq
\sup_{a \in \dR} \mathrm{Ent}_{\mu_2}((f+a)^2)
\leq
\frac{5}{2} \NRM{(f-\mu_2(f))^2}_{\Phi_1^1} .
$$
Thus, the previous result gives that for $t \geq 4Ce$,
$$
\mathrm{Ent}_{\mu_2}((\PT{t}f)^2) \leq
\frac{15t}{16C} e^{-\frac{t}{4eC}}
\sup_{a \in \dR} \mathrm{Ent}_{\mu_2}((f+a)^2) ,
$$
where $C$ is the logarithmic Sobolev constant of $\mu_2$.
Now, using the Rothaus inequality (see \cite{rothaus})
$\sup_{a \in \dR} \mathrm{Ent}_{\mu_2}((f+a)^2) \leq
\mathrm{Ent}_{\mu_2}\left((f-\mu_2 (f))^2\right) +
2 \mu_2\left((f-\mu_2 (f))^2\right)$,
we have
\begin{equation} \label{eq:ent}
\mathrm{Ent}_{\mu_2}((\PT{t}f)^2)
\leq
\frac{15t}{16C} e^{-\frac{t}{4eC}}
\left(\mathrm{Ent}_{\mu_2}\left((f-\mu_2 (f))^2\right)+
2 \mu_2\left((f-\mu_2 (f))^2\right)\right)
\end{equation}
which can be improved for $f\geq 0$ using Kulback's inequality
$\var_{\mu_2}(f) \leq \mathrm{Ent}_{\mu_2}(f^2)$.
As far as we know the bound (\ref{eq:ent}) was not known. Indeed,
the logarithmic
Sobolev inequality is usually used in case of diffusion semi-group
(see \textit{e.g.} \cite{ane})
to prove exponential decay of entropy, \textit{i.e.} that for any $t$,
$$
\mathrm{Ent}_{\mu}(\PT{t}f)
\leq e^{-\frac{t}{C}} \mathrm{Ent}_\mu(f) .
$$
On the other hand there does not exist any constant $k<\infty$ such that
for any function $f$,
$\sup_{a \in \dR} \mathrm{Ent}_\mu((f+a)^2) \leq k \mathrm{Ent}_\mu(f^2)$,
or equivalently
$\NRM{(f-\mu(f))^2}_{\Phi_1^1} \leq k \mathrm{Ent}_\mu(f^2)$.
Indeed, on the space $\{0,1\}$ with the symmetric Bernoulli measure, consider
the function $f(0)=-1$ and $f(1)=1$ for which $(f-\mu(f))^2 \equiv \ind$ and
$\mathrm{Ent}_\mu(f^2) = 0$.

Thus we will consider the functional
$$
A(f) \equiv \mathrm{Ent}_{\mu_2}(f^2) + \mu_2 \left(f-\mu_2 (f)\right)^2 .
$$
Then the bound (\ref{eq:ent}), for all $t>T$ with some $T\in (0,\infty)$,
can be written as follows
$$
A(\PT{t}f) \leq e^{-mt} A(f)
$$
with some $m\in(0,\infty)$.
With $\omega\in (0, m)$, define
$$
\mathcal{A}_\omega(f) \equiv \sup_{s\in[0,T]} A(\PT{s}f)e^{\omega s}
$$
and for $\omega\in[0,m]$ define
$$
\mathcal{B}_\omega(f) \equiv \frac 1T \int_0^T A(\PT{s}f)e^{\omega s}ds .
$$

\begin{proposition}
Suppose, with some $m, T\in(0,\infty)$, for all $t\geq T$, one has
$$
A(\PT{t}f) \leq e^{-mt} A(f) .
$$
Then the functionals
$\mathcal{A}_\omega$ and $\mathcal{B}_\omega$  are exponentially decaying,
that is for any $t\geq 0$
$$
\mathcal{A}_\omega(\PT{t}f) \leq e^{-\omega t}\mathcal{A}_\omega(f)
$$
and
$$
\mathcal{B}_\omega(\PT{t}f) \leq e^{-\omega t}\mathcal{B}_\omega(f) .
$$
\end{proposition}

\begin{proof}
If $t\geq T$, the statements are clear.
For $\mathcal{A}_\omega$ and $0\leq t\leq T$ note that
\begin{eqnarray*}
\mathcal{A}_\omega(\PT{t}f)
& = &
\sup_{s\in[0,T]} A(\PT{s+t}f)e^{\omega s}
=
e^{-\omega t}\sup_{s\in[t,T+t]} A(\PT{s}f)e^{\omega s} \\
& = &
e^{-\omega t}\max \big(\sup_{s\in[t,T]} A(\PT{s}f)e^{\omega s},
\sup_{s\in[T,T+t]} A(\PT{s}f)e^{\omega s}\big) .
\end{eqnarray*}
Since for $s\in[T,T+t]$
$$
A(\PT{s}f)e^{\omega s}
\leq
e^{-m T} A(\PT{s-T}f)
\leq
e^{-m T+\omega T} \big(A(\PT{s-T}f)e^{\omega (s-T)}
\big)
$$
we get for $t\in[0,T]$ and $\omega\leq m$
$$
\sup_{s\in[T,T+t]} A(\PT{s}f)e^{\omega s}
\leq
\sup_{s\in[0,T]} A(\PT{s}f)e^{\omega s} .
$$
This together with the previous considerations concludes the arguments
for exponential decay of the first functional.
In case of $\mathcal{B}_\omega$, for $0\leq t\leq T$, we have
$$
T\mathcal{B}_\omega(\PT{t}f) \equiv \int_0^T A(\PT{s+t}f)e^{\omega s}ds =
e^{-\omega t}\int_t^{T+t} A(\PT{s}f)e^{\omega s}ds .
$$
Next we note that
$$
\int_t^{T+t} A(\PT{s}f)e^{\omega s}ds = \int_t^{T} A(\PT{s}f)e^{\omega s}ds
+\int_T^{T+t} A(\PT{s}f)e^{\omega s}ds .
$$
To complete the proof it is sufficient to note that
\begin{eqnarray*}
\int_T^{T+t} A(\PT{s}f)e^{\omega s}ds
& \leq &
e^{-mT+\omega T}\int_T^{T+t} A(\PT{s-T}f)e^{\omega (s-T)}ds \\
& = &
e^{-mT+\omega T}\int_0^{t} A(\PT{s}f)e^{\omega (s)}ds \\
& \leq &
\int_0^{t} A(\PT{s}f)e^{\omega (s)}ds .
\end{eqnarray*}
\end{proof}

In particular we have thus shown that if an (a priori non-convex)
functional decays monotonously exponentially fast for large times,
then by averaging over "a characteristic time of relaxation" we can get a
globally monotone functional.


\section{Orlicz-Sobolev and $\Phi$-Sobolev ine\-qua\-lities}

In this section we provide a link between the Orlicz-Sobolev inequality
and the $\Phi$-Entropy bound introduced by Chafa{\"\i} \cite{chafai}
and the \textit{additive} $\Phi$-Sobolev
inequality studied in \cite{barthe-cattiaux-roberto}.

Given a closed interval $\cal I$ of $\dR$ and
a convex function $\Phi : {\cal I} \rightarrow \dR$,
a probability measure $\mu$ on $\dR^n$ satisfies a
$\Phi$-Sobolev inequality if there exists a constant $C_\Phi$
such that for every smooth function $f : \dR^n \rightarrow {\cal I}$,
\begin{equation}\label{eq:chafai}
\tag*{(\textbf{Ent$^\Phi$-S})}
\mathrm{Ent}_\mu^\Phi (f)
\leq C_\Phi \int \Phi''(f) |\nabla f|^2 d\mu
\end{equation}
where
$$
\mathrm{Ent}_\mu^\Phi (f) : = \int \Phi(f) d\mu - \Phi \left( \int f d\mu \right)
$$
In \cite{chafai}, it is proved that such an inequality is equivalent to the
exponential decay of $\mathrm{Ent}_\mu^\Phi ( \PT{t} f)$.

On the other hand, given a non-decreasing function
$\varphi : (0,+\infty) \rightarrow \dR$ continuously differentiable, we define
$\Phi(x)=x \varphi(x)$ and we assume that $\Phi$ can be extended to $0$
and is convex. A probability measure $\mu$ on $\dR^n$ satisfies an additive
$\Phi$-Sobolev inequality if there exists a constant $C_\Phi$
such that for every smooth function $f : \dR^n \rightarrow \dR$,
\begin{equation}\label{eq:additive-chafai}
\tag*{({\textbf{$\Phi$-S}})}
\int \Phi(f^2) d\mu - \Phi \left( \int f^2 d\mu \right)
\leq C_\Phi \int |\nabla f|^2 d\mu .
\end{equation}

We start with the following general fact.

\begin{proposition} \label{prop:poincare}
Let $\Phi(x)=x \varphi(x)$ be a $\C{2}$ Young function, with
$\varphi : (0,+\infty) \rightarrow \dR$ non-decreasing.
Assume that the probability measure $\mu$
on $\dR^n$ satisfies for any smooth function $f$,
$$
\int \Phi(f^2) d\mu - \Phi \left( \int f^2 d\mu \right)
\leq C_\Phi \int |\nabla f|^2 d\mu ,
$$
for some constant $C_\Phi$ independent of $f$. Then, for any smooth
function $g$, for any $a >0$,
$$
\Phi''(a) \var_\mu(g) \leq \frac{C_\Phi}{2a} \int |\nabla g|^2 d\mu .
$$
In particular, if $\Phi'' \neq 0$,
$\mu$ satisfies a Poincar{\'e} inequality with constant
$C_p \leq \inf_{a > 0} \frac{C_\Phi}{2a \Phi''(a)}$.
\end{proposition}

The previous result states that the Poincar{\'e} inequality holds as
far as the additive $\Phi$-Sobolev inequality holds and $\Phi'' \neq 0$.

\begin{proof}
Given a smooth non negative function $f$ on $\dR^n$,
the additive $\Phi$-Sobolev inequality applied to $\sqrt f$ leads to
$$
\int \Phi(f) d\mu - \Phi \left( \int f d\mu \right)
\leq \frac{C_\Phi}{4} \int \frac{|\nabla f|^2}{f} d\mu .
$$
Now, given a smooth bounded function $g$ with $\mu(g)=0$ and $a >0$,
$a + \varepsilon g \geq 0$ for $\varepsilon$ small enough.
Then the previous inequality applied to $a + \varepsilon g$ and a Taylor
expansion at the second order for $\Phi$ gives the result when
$\varepsilon$ tends to $0$.
\end{proof}

\begin{remark}
In \cite[section 1.2]{chafai}, the same result is proved for the
$\Phi$-Entropy bound ({\textbf{Ent$^\Phi$-S}}).

On the other hand, $\Phi''  \equiv 0$ is equivalent to
$\varphi(x)=a - (b/x)$, $(a,b) \in \dR \times \dR^+$. In that case the additive
$\Phi$-Sobolev inequality is trivial.
\end{remark}

Now we give a link between the modified Orlicz-Sobolev inequality
\eqref{in:mos} and the $\Phi$-Entropy bound ({\textbf{Ent$^\Phi$-S}}).

\begin{proposition} \label{prop:c2}
Let $\Phi$ be a Young function. Assume that the probability measure $\mu$
on $\dR^n$ satisfies a $\Phi$-Entropy bound ({\textbf{Ent$^\Phi$-S}})
with constant $C_\Phi$. Then, it satisfies a modified Orlicz-Sobolev
inequality \eqref{in:mos} with the same constant  $C_\Phi$.
\end{proposition}

\begin{proof}
For every smooth function $f : \dR^n \rightarrow \dR$
apply the $\Phi$-Entropy bound ({\textbf{Ent$^\Phi$-S}})
to
$(f-\mu(f))/\NRM{f-\mu(f)}_{\Phi}$ to get
\begin{eqnarray*}
&&
\int \Phi \left( \frac{f-\mu(f)}{\NRM{f-\mu(f)}_{\Phi}} \right) d\mu
- \Phi \left( \int  \frac{f-\mu(f)}{\NRM{f-\mu(f)}_{\Phi}} d\mu \right)\\
&&  \qquad \qquad \qquad \qquad \leq C_\Phi
\int \Phi''\left(\frac{f-\mu(f)}{\NRM{f-\mu(f)}_{\Phi}} \right)
\frac{|\nabla f|^2}{\NRM{f-\mu(f)}_{\Phi}^2} d\mu .
\end{eqnarray*}
Since $\Phi(0)=0$ and
$\int \Phi \left( \frac{f-\mu(f)}{\NRM{f-\mu(f)}_{\Phi}} \right) d\mu=1$,
we get the expected result.
\end{proof}

\begin{remark}
As a consequence of this result and using Theorem \ref{th:decay}, we get that
if $\mathrm{Ent}_\mu^\Phi ( \PT{t} f)$ decays
exponentially fast in time, then $\NRM{\PT{t}f-\mu(f)}_{\Phi}$
decays exponentially fast.
\end{remark}

Next we give a similar result involving the additive $\Phi$-Sobolev inequality
({\textbf{$\Phi$-S}})
and the Orlicz-Sobolev inequality ({\textbf{O-S}}).
Note that for a Young function $\Phi$,
the assumption $\Phi(x)/x \nearrow \infty$ when $x$ goes to
infinity and $\Phi'(0)>0$ insure that the equation $x \Phi'(x)=1$ has a
unique solution, see \cite[section 2.4]{rao-ren}.

\begin{proposition} \label{prop:phi}
Let $\Phi$ be a $\C{2}$ Young function with $\Phi'(0)>0$.
Assume that $\Phi(x)=x \varphi (x)$ for a non-decreasing function
$\varphi$ defined on $(0,\infty)$ and such that
$\lim_{+\infty} \varphi = + \infty$.
Denote by $k_0$ be the unique solution of $k_0 \Phi'(k_0)=1$.
Let $\mu$ be a probability measure on $\dR^n$.
Assume that there exists a constant $C_\Phi$
such that for every smooth function $f: \dR^n \rightarrow \dR$,
$$
\int \Phi(f^2) d\mu - \Phi \left( \int f^2 d\mu \right)
\leq C_\Phi \int |\nabla f|^2 d\mu .
$$
Then, for any smooth function $f: \dR^n \rightarrow \dR$, for any $a>0$,
$$
\NRM{(f-\mu(f))^2}_{\Phi} \leq
\frac{C_\Phi}{k_0} \left( \frac{1}{2 a \Phi''(a)} + \frac{1}{\Phi'(0)} \right)
\int |\nabla f|^2 d\mu .
$$
\end{proposition}

\begin{remark}
Note that since $\lim_{x\to +\infty} \varphi = + \infty$, there
exists $a>0$ such that $\Phi''(a) >0$.

On the other hand,
it is easy to get rid of the assumption $\Phi'(0)>0$. Indeed, assume that
$\Phi'(0)=0$ and defined $\Phi_\lambda(x)=\Phi(x)+ \lambda |x|$ for
$\lambda >0$. Then, on one hand $\Phi_\lambda'(0)=\lambda >0$.
On the other hand,
if an additive $\Phi$-Sobolev inequality holds, then a
$\Phi_\lambda$-Sobolev inequality holds,
with the same constant. So the previous
Proposition applies to $\Phi_\lambda$:
for any smooth function $f$, for any $a>0$,
$$
\NRM{(f-\mu(f))^2}_{\Phi_\lambda} \leq
\frac{C_\Phi}{k_0(\lambda)}
\left( \frac{1}{2a \Phi''(a)} + \frac{1}{\lambda} \right)
\int |\nabla f|^2 d\mu .
$$
Since $\Phi \leq \Phi_\lambda$,
$\NRM{(f-\mu(f))^2}_{\Phi} \leq \NRM{(f-\mu(f))^2}_{\Phi_\lambda}$. This leads
to
$$
\NRM{(f-\mu(f))^2}_{\Phi} \leq
\frac{C_\Phi}{k_0(\lambda)}
\left( \frac{1}{2 a \Phi''(a)} + \frac{1}{\lambda} \right)
\int |\nabla f|^2 d\mu  ,
$$
for any  $\lambda >0$, any $a>0$ and any function $f$.
Note that $k_0(\lambda) \rightarrow 0$ when $\lambda$ tends to $\infty$.
\end{remark}

\begin{proof}
Let $\widetilde \Phi(x):=\Phi(k_0 x)$, so its complementary
function is $(\widetilde \Phi)^* (x) = \Phi^* (x/k_0)$ where
$\Phi^*$ is the complementary function of $\Phi$. Now
$(\widetilde \Phi, (\widetilde \Phi)^*)$ is a normalized complementary pair
of Young functions.
Following \cite{rao-ren} define the following modified Luxembourg norm
$$
|| f ||_{\widetilde \Phi} =
\inf \{ \lambda ; \int \widetilde \Phi \left( \frac{f}{\lambda} \right) d\mu
\leq \widetilde \Phi(1) \} .
$$
Note that $|| \ind ||_{\widetilde \Phi}=1$.
By \cite[Proposition 1 in section 3.3]{rao-ren}, we know that
\begin{equation}\label{eq:holder2}
\int |g| d \mu
\leq
|| g ||_{\widetilde \Phi}
|| \ind ||_{(\widetilde \Phi)^*}
=
|| g ||_{\widetilde \Phi}
\qquad \qquad \forall g \in \L_{\widetilde \Phi}(\mu) .
\end{equation}
It is important to introduce this modified norm in order to have the latter
inequality with a factor $1$ in front of the r.h.s and not $2$ as in
the standard inequality \eqref{in:holderp}.

Now let $f: \dR^n \rightarrow \dR$ be a smooth function.
From the additive $\phi$-Sobolev inequality applied to
$\sqrt{k_0} (f- \mu(f)) / || (f- \mu(f))^2 ||_{\widetilde \Phi}^{1/2}$,
we get
\begin{eqnarray*}
&& \int \widetilde \Phi \left(
\frac{(f- \mu(f))^2}{|| (f- \mu(f))^2 ||_{\widetilde \Phi}}  \right) d\mu
-
\widetilde \Phi
\left( \int \frac{(f- \mu(f))^2}{|| (f- \mu(f))^2 ||_{\widetilde \Phi}} d\mu
\right) \\
&& \qquad \qquad \qquad \qquad \qquad \qquad \qquad \qquad
\leq C_\Phi k_0
\int \frac{|\nabla f|^2}{|| (f- \mu(f))^2 ||_{\widetilde \Phi}}d\mu.
\end{eqnarray*}
Since $\int \widetilde \Phi \left(
\frac{(f- \mu(f))^2}{|| (f- \mu(f))^2 ||_{\widetilde \Phi}} \right) d\mu =\widetilde \Phi(1)$,
it follows that
$$
\widetilde \Phi (1)
-\widetilde \Phi \left(\frac{\int g d\mu }{N_{\widetilde \Phi}(g)} \right)
\leq
C_\Phi k_0
\int \frac{|\nabla g|^2}{|| g ||_{\widetilde \Phi}}d\mu,
$$
where $g:=(f-\mu(f))^2$.
A Taylor expansion of $\widetilde \Phi$ up to the second order, between
$1$ and $\frac{\int g d\mu }{|| g ||_{\widetilde \Phi}}$,
and convexity of $\widetilde \Phi$, give that
\begin{eqnarray*}
\widetilde \Phi (1)
-\widetilde \Phi \left(\frac{\int g d\mu }{|| g ||_{\widetilde \Phi}} \right)
& = &
\left(1-\frac{\int g d\mu }{|| g ||_{\widetilde \Phi}} \right)
\widetilde \Phi' \left(\frac{\int g d\mu }{|| g ||_{\widetilde \Phi}} \right)
\!+\! \frac 12 \left(1-\frac{\int g d\mu }{|| g ||_{\widetilde \Phi}} \right)^2 \!\!
\widetilde \Phi''(\theta) \\
& \geq &
\left(1-\frac{\int g d\mu }{|| g ||_{\widetilde \Phi}} \right)
\widetilde \Phi' \left(\frac{\int g d\mu }{|| g ||_{\widetilde \Phi}}\right) \\
& \geq &
\left(1-\frac{\int g d\mu }{|| g ||_{\widetilde \Phi}} \right)
\widetilde \Phi'(0) ,
\end{eqnarray*}
where $\theta \in (0,1)$ (recall that from \eqref{eq:holder2},
$\frac{\int |g| d\mu }{|| g ||_{\widetilde \Phi}} \leq 1$).
This leads to
$$
{|| (f-\mu(f))^2 ||_{\widetilde \Phi}} \leq
\frac{C_\Phi k_0 }{\widetilde \Phi'(0)} \int |\nabla f|^2 d\mu
+ \var_\mu(f).
$$
Since $\lim_{x\to+\infty} \varphi(x) = + \infty$, there
exists $a>0$ such that $\Phi''(a) >0$. Choose such an $a$.
>From Proposition \ref{prop:poincare}, $\mu$ satisfies a Poincar{\'e} inequality
with constant less than $C_\Phi/(2a \Phi''(a))$.
On the other hand
$|| (f-\mu(f))^2 ||_{\widetilde \Phi} = k_0
\NRM{(f-\mu(f))^2}_{\Phi/ \widetilde \Phi(1)}$
and $\widetilde \Phi'(0)=k_0 \Phi'(0)$. Thus, for any smooth
function $f : \dR^n \rightarrow \dR$,
$$
\NRM{(f-\mu(f))^2}_{\Phi/ \widetilde \Phi(1)}
\leq \frac{C_\Phi}{k_0}  \left( \frac{1}{2 a \Phi''(a)} + \frac{1}{\Phi'(0)} \right)
\int |\nabla f|^2 d\mu.
$$
The result follows from the fact that
$\NRM{(f-\mu(f))^2}_{\Phi} \leq \NRM{(f-\mu(f))^2}_{\Phi/\widetilde \Phi(1)}$
since $\Phi \leq \Phi/\widetilde \Phi(1) $ (recall that
$\widetilde \Phi(1) + (\widetilde \Phi)^*(1)= 1$).

\noindent For all $a>0$ such that $\Phi''(a)=0$, the result is trivial.
This ends the proof.
\end{proof}

\noindent Proposition \ref{prop:phi} allows us to give a criterium for the
$\Phi$-Sobolev inequality to hold. This completes
\cite[Theorem 26]{barthe-cattiaux-roberto}.

\begin{theorem}\label{th:bcrphi}
Let $\Phi(x)=x \varphi(x)$ be a $\C{2}$ Young function
with $\varphi$ non decreasing, concave, with $\varphi(0)>0$
and such that $\lim_{+\infty} \varphi = + \infty$. Denote by  $k_0$
the unique solution of $k_0 \Phi'(k_0)=1$. Assume that there exist constants
$\gamma, \kappa$ and such that for all $x,y >0$ one has
$$
x \varphi'(x) \leq \gamma \qquad \mbox{and} \qquad
\varphi(xy) \leq \kappa + \varphi(x) + \varphi(y),
$$
and a constant $\lambda \geq 2$ such that for every $x \geq 2\lambda$, one has
$\lambda \varphi(x/\varphi(x)) \geq \varphi(x)$.

Let $\mu$ be a probability measure on $\dR^n$ satisfying the Poincar\'e inequality
with constant $C_P$ and $C_\Phi$  the optimal constant
such that
for every smooth function $f: \dR^n \rightarrow \dR$ one has
\begin{equation} \label{eq:as}
\int \Phi(f^2) d\mu - \Phi \left( \int f^2 d\mu \right)
\leq C_\Phi \int |\nabla f|^2 d\mu .
\end{equation}
Then, for any $a>0$,
$$
\frac{k_0a \Phi''(a) \varphi(0)}{8\lambda(\varphi(0) + 2a \Phi''(a))}
\widetilde B(\Phi) \leq C_\Phi \leq
(18\gamma C_p + 24(1+\frac{M}{\varphi(8)})) \widetilde B(\Phi)
$$
where $\widetilde B(\Phi)$ is the smallest constant so that for every
$A \subset \dR^n$ with $\mu(A) < \frac 12$
$$
\mu(A) \varphi \left( \frac{2}{\mu(A)} \right)
\leq \widetilde B(\Phi) \mathrm{Cap}_\mu (A) .
$$
\end{theorem}

\begin{proof}
The upper bound on $C_\Phi$
follows from \cite[Theorem 26]{barthe-cattiaux-roberto}.
\newline
Assume that the additive $\Phi$-Sobolev inequality
\eqref{eq:as} holds. Then, by Proposition \ref{prop:phi},
for every smooth function $f : \dR^n \rightarrow \dR$, every $a>0$,
$$
\NRM{(f-\mu(f))^2}_{\Phi}
\leq \frac{C_\Phi}{k_0} \left( \frac{1}{2 a \Phi''(a)} + \frac{1}{\varphi(0)}\right)
\int |\nabla f|^2 d\mu.
$$
Then, by Theorem \ref{th:cos} we get,
\begin{equation}\label{eq:step1}
\frac 18 B( \Phi) \leq
\frac{C_\Phi}{k_0} \left(\frac{1}{2 a \Phi''(a)} + \frac{1}{\varphi(0)} \right)
\end{equation}
where $B(\Phi)$ is
the smallest constant so that for every
$A \subset \dR^n$ with $\mu(A) < \frac 12$
$$
\frac{1}{\Phi^{-1}\left( \frac{1}{\mu(A)} \right)} = \NRM{\ind_A}_{\Phi}
\leq B(\Phi) \mathrm{Cap}_\mu (A) .
$$
By our assumption on $\varphi$,
$\Phi \left( \frac{x}{\varphi(x)} \right) =
x \frac{\varphi(x/\varphi(x))}{\varphi(x)} \geq \frac 1\lambda x$
for all $x \geq 2 \lambda$. Thus, since $\lambda \geq 2$ and $\varphi$
is non-decreasing, for all $y \geq 2$
$$
\Phi^{-1}(y) \leq \frac{\lambda y}{\varphi(\lambda y)}
\leq \lambda \frac{y}{\varphi(2y)} .
$$
It follows that $\widetilde B(\Phi) \leq \lambda  B(\Phi)$. This together
with \eqref{eq:step1} achieves the proof.
\end{proof}

\subsection*{$\Phi$-\textbf{S} and \textbf{O-S}
Inequalities in Infinite Dimensions}

It is not difficult to check that
$\Phi(x)=|x|\big(\log(\eta + |x|)\big)^\beta$, $\beta \in(0,1]$, $\eta >1$,
satisfies the hypothesis of Theorem \ref{th:bcrphi}.

Following a remark of \cite{barthe-cattiaux-roberto} we note that
$$
\mu \left(\Phi(f^2)\right) - \Phi(\mu (f^2))
=
\inf_{t>0} \mu\big(\Phi(f^2) -\Phi(t)-\Phi'(t)(\mu(f^2) -t)\big)
$$
By convexity of $\Phi$ one has
$\Phi(f^2) -\Phi(t)-\Phi'(t)(\mu(f^2) -t) \geq 0$
which implies the following Mild Perturbation Property (MPP) for additive
$\Phi$-Sobolev inequality.

\begin{proposition}
Let $d\mu =\rho d\nu$ with
$\delta U\equiv \sup (\log \rho) - \inf(\log\rho) < \infty$ and assume that
$$
\int \Phi(f^2) d\nu - \Phi\left( \int f^2d \nu \right)
\leq C \int |\nabla f|^2 d\nu .
$$
Then
$$
\int \Phi(f^2) d\mu - \Phi\left( \int f^2d \mu \right)
\leq C e^{\delta U} \int |\nabla f|^2 d\mu  .
$$
\end{proposition}

The additive $\Phi$-Sobolev inequality, with the $\Phi$ as described above,
was in particular established for products of
$\mu_\alpha$ measures with suitable $\alpha\in(1,2)$.
Using MPP one can construct a compatible family of
finite dimensional expectations $E_\Lambda$
(with partially ordered indices $\Lambda$)
for which additive $\Phi$-Sobolev inequality also holds.
By definition for the corresponding Gibbs measure
${\boldsymbol\mu} (E_\Lambda) = {\boldsymbol\mu} $
and one has the following simple conditioning property
\begin{eqnarray*}
{\boldsymbol\mu} \left(\Phi(f^2)\right)- \Phi({\boldsymbol\mu} (f^2))
&=&
{\boldsymbol\mu} \left[E_\Lambda\left(\Phi(f^2)\right) -
\Phi\left(E_\Lambda (f^2)\right)\right] \\
&&\quad +
{\boldsymbol\mu} \left(\Phi\left[E_\Lambda (f^2)\right] \right)-
\Phi\left(\mu\left[E_\Lambda (f^2)\right]\right) .
\end{eqnarray*}
With these two facts in mind, under suitable mixing condition,
one can follow closely
the strategy originally invented for the proof of Logarithmic
Sobolev Inequality  (\textit{cf.} \cite{GZ03}) to proof the following result

\begin{theorem}
Suppose a local specification is mixing and satisfies
$\Phi$-Sobolev inequality. Then the unique
Gibbs measure ${\boldsymbol\mu}$ satisfies
$$
{\boldsymbol\mu} \left(\Phi(f^2)\right)
- \Phi({\boldsymbol\mu} (f^2))
\leq C
{\boldsymbol\mu}\left(|\nabla f|^2 \right)
$$
with a constant $C$ independent of a function $f$.
\end{theorem}

This provides a large family of  nontrivial examples of (non-product)
measures on infinite dimensional spaces satisfying additive
$\Phi$-Sobolev inequality.

We remark that by inserting into such the inequality a function $f/||f||_2$ and setting $F(x) \equiv  \big(\log(\eta + |x|)\big)^\beta -
 \big(\log(\eta + 1)\big)^\beta $, we arrive at the following
$F$-Sobolev inequality,
\begin{corollary}
\begin{equation}
\tag*{({\textbf{F-S}})}
\int f^2F \left(\frac {f^2}{ {\boldsymbol\mu} (f^2)} \right) d{\boldsymbol\mu}
\leq
C {\boldsymbol\mu} \int |\nabla f|^2  d{\boldsymbol\mu}
\end{equation}
for the Gibbs measure ${\boldsymbol\mu}$.
\end{corollary}

Finally we note that by the same arguments as the ones used to prove
Proposition \ref{prop:phi},
we get the following Orlicz-Sobolev inequality for infinite dimensional
Gibbs measures
\begin{corollary}
$$
\NRM{(f-{\boldsymbol\mu} (f))^2}_{\Phi}
\leq
c {\boldsymbol\mu} \left(|\nabla f|^2 \right)
$$
with a constant $c$ independent of a function $f$.
\end{corollary}


\section{Orlicz-Sobolev and Nash-type inequalities}

In this section we prove that the
Orlicz-Sobolev inequality is equivalent,
up to some constants, to a Nash-type inequality.
This give
new results on the decay to equilibrium of the
semi-group (see next section).

\begin{theorem}\label{th:nash}
Let $\Phi$ and $\Psi(x)=\frac{x^2}{\psi(|x|)}$ be two $N$-functions
with $\psi: \dR^+ \to \dR^+$ increasing, satisfying $\psi(0)=0$ and
$\lim_{+\infty}\psi =+ \infty$. Assume that the probability measure
$\mu$ on $\dR^n$ satisfies for any smooth function $f: \dR^n
\rightarrow \dR$
$$
\NRM{(f-\mu(f))^2}_\Phi \leq C_\Phi \int |\nabla f|^2  d \mu .
$$
Then, for any function $f$,
\begin{equation} \label{eq:nash}
\var_\mu(f) \theta \left( \frac12
\frac{\var_\mu(f)}{\NRM{f-\mu(f)}_\Psi^2} \right) \leq  4 C_\Phi
\int |\nabla f|^2  d \mu .
\end{equation}
where $\theta={\Phi^*}^{-1} \circ \Psi \circ \psi^{-1}$ (here
$\Phi^*$ is the complementary pair of $\Phi$;
${\Phi^*}^{-1}$ and $\psi^{-1}$ stand for the inverse function of
$\Phi^*$ and $\psi$ respectively).
\end{theorem}

\begin{remark}
Note that by our assumption on $\psi$, $\psi^{-1}$ is well defined
on $\dR_+$ onto $\dR_+$.

Furthermore, in order to deal with explicit functions, one can easily
see that under the assumption of the Theorem,
$\Psi(x) \leq \frac{1}{\psi(1)}(x+x^2)$ in such a way that
$\var_\mu(f)/\NRM{f-\mu(f)}_\Psi^2 \geq c$ for some constant $c$ (see
Lemma \ref{lemma_injection_L1}). Thus one has only to consider the behavior of
$\theta$ (or equivalently to $\Phi$, $\Psi$ and $\psi$) away from $0$.
\end{remark}

\begin{remark}
We will call the inequality
$$
\var_\mu(f) \theta \left( \frac12
\frac{\var_\mu(f)}{\NRM{f-\mu(f)}_\Psi^2} \right) \leq  C \int
|\nabla f|^2  d \mu ,
$$
a Nash-type inequality since for $\Phi(x)=|x|^{(d)/(d-2)}$,
$\Psi(x)=\psi(x)=x$ (and thus $\theta(x)=c_d |x|^{\frac 2d}$ for
some constant $c_d$), it reads for any $f$ with $\mu(f)=0$ as
$$
\NRM{f}_2^{1+\frac 2d} \leq C' \NRM{\nabla f}_2
\NRM{f}_{1}^{\frac 2d}
$$
which is the standard Nash inequality (\cite{nash}).
\end{remark}

\begin{proof}
The proof is a generalization of \cite[Proposition
10.3]{bakry-coulhon-ledoux-sc}, see also \cite{roberto}. Let $f$ be
a function with $\mu(f)=0$ and $\NRM{f}_\Psi=1$ in such a way that
$\int \Psi(f)d \mu = 1$. Fix a parameter $t >0$. Denote by $\Phi^*$
the complementary function of $\Phi$. From \eqref{in:holderp}, if
$(f,g)  \in \L_\Phi \times \L_{\Phi^*}$, $\int |fg| d\mu \leq 2
\NRM{f}_\Phi \NRM{f}_{\Phi^*}$. Hence,
\begin{eqnarray*}
\var_\mu(f)
& = &
\int f^2 \ind_{|f| <t} d\mu + \int f^2 \ind_{|f| \geq t} d\mu \\
&\leq & \int \Psi(f) \psi(|f|) \ind_{|f| <t} d\mu +
2  \NRM{f^2}_\Phi \NRM{\ind_{|f| \geq t}}_{\Phi^*} \\
&\leq & \psi(t) \int \Psi(f) d\mu + \frac{2C_\Phi}{{\Phi^*}^{-1}
(1/\mu(|f| \geq t))}
 \int |\nabla f|^2  d \mu .
\end{eqnarray*}
Now by Chebychev inequality (recall that $\Phi$ is an even function) we have
$$
\mu(|f| \geq t) = \mu (\Psi(f) \geq \Psi(t)) \leq \frac{1}{\Psi(t)}
\int \Psi(f) d\mu = \frac{1}{\Psi(t)}.
$$
It follows that
for any $t>0$,
$$
\var_\mu(f) \leq \psi(t) + \frac{2 C_\Phi}{{\Phi^*}^{-1}(\Psi(t))}
\int |\nabla f|^2 d \mu .
$$
Now choose $t$ such that $\psi (t)= \frac 12 \var_\mu(f)$. We get
$$
{\Phi^*}^{-1}(\Psi(t)) \var_\mu(f)
\leq 4C_\Phi \int |\nabla f|^2
d\mu .
$$
This gives the expected result by homogeneity.
\end{proof}

\begin{example} \label{ex:nash}
Let $\alpha \in [1,2]$,
$\beta = 2(1-\frac1\alpha) \in [0,1]$
and define  the probability measure on $\dR^n$:
$d\mu_\alpha^n(x) =
Z_\alpha^{-n} \exp \{-\sum_{i=1}^n |x_i|^\alpha \}dx$.
For any $\gamma \geq 1$ define
$\Phi_{\beta}^\gamma (x) = |x| (\log(\gamma+|x|)^\beta$ and
$\Phi_{\beta,2}^\gamma (x) = \Phi_\beta^\gamma (x^2)$.
From Corollary \ref{cor:alpha}, point $(i)$ of
Lemma \ref{lem:eq} and the general fact
that $\NRME{f^2}{\Phi_\beta^\gamma}=\NRM{f}^2_{\Phi_{\beta,2}^\gamma}$,
there exists a constant $C$ (independent of $n$) such that
for any smooth function $f : \dR^n \rightarrow \dR$,
$$
\NRM{(f-\mu_\alpha^n(f))^2}_{\Phi_{\beta}^\gamma}
\leq
C (1+e (\log \gamma)^\beta) \int |\nabla f|^2 d \mu_\alpha^n.
$$
Using similar computation than in the proof of Inequality \eqref{eq:inverse},
it is not difficult to see that for any $\varepsilon >0$, there exists a constant
$C_\varepsilon$ (depending also on $\beta$ and $\gamma$) such that
for any $x \geq \varepsilon$,
$$
C_\varepsilon^{-1} \log(1+x)^\beta
\leq
{{\Phi_{\beta}^\gamma}^*}^{-1} (x)
\leq
C_\varepsilon \log(1+x)^\beta .
$$
Now define for $x \geq 0$ and $\delta \in (0,1)$, $\psi(x)=\left(\log(1+x)\right)^\delta$.
One can easily see that $\Psi(x):=x^2/\psi(x)$ is a $N$-function.
We deduce that there exists $C_\varepsilon' >0$ such that for any $x \geq \varepsilon$,
$$
{C_\varepsilon'}^{-1} x^\frac \beta \delta
\leq
\theta (x)
\leq
C_\varepsilon' x^\frac \beta \delta
$$
where $\theta : = {{\Phi_{\beta}^\gamma}^*}^{-1} \circ \Psi \circ \psi^{-1}$.
Theorem \ref{th:nash} implies that
there exists a constant $C'$ (independent on $n$ and possibly depending on
$\beta$, $\delta$) such that for any function $f : \dR^n \rightarrow \dR$,
$$
\var_{\mu_\alpha^n}(f)^{1+\frac{\beta}{\delta}}
\leq C'
\NRM{f-\mu_\alpha^n(f)}_{\Psi}^\frac{2 \beta}{\delta}
\int |\nabla f|^2  d \mu_\alpha^n .
$$

\bigskip

If we choose instead $\widetilde \psi(x)= e^{\left(\log(1+x)\right)^\delta}-1$ for
$\delta \in (0,1)$,
$\widetilde \Psi(x)=x^2/\widetilde \psi(x)$
is again a $N$-function. It follows in this case that there exists a constant
$C''_\varepsilon>0$ such that for any
$x \geq \varepsilon$,
$$
{C_\varepsilon'}^{-1} \log(1+x)^\frac \beta \delta
\leq
\widetilde \theta (x)
\leq
C_\varepsilon' \log(1+x)^\frac \beta \delta
$$
where $\widetilde \theta : = {{\Phi_{\beta}^\gamma}^*}^{-1}
\circ \widetilde \Psi \circ {\widetilde \psi}^{-1}$.
In turn, Theorem \ref{th:nash} implies that
$$
\var_{\mu_\alpha^n}(f) \log \left( 1+ \frac 12
\frac{\var_{\mu_\alpha^n}(f)}
{\NRM{f-\mu_\alpha^n(f)}_{\Psi}^2} \right)^\frac{\beta}{\delta}
\leq C''
\int |\nabla f|^2  d \mu_\alpha^n
$$
for some constant $C''$ independent on $n$ and $f$.
\end{example}

It is natural to ask for the equivalence between the Orlicz-Sobolev
inequality and the Nash-type inequality in Theorem \ref{th:nash}.
It seems (almost for us) to be difficult to prove
directly this equivalence. However  it is possible to achieve
that with the help of an intermediate inequality as follows.

As a first step, we consider the following equivalent form of the Nash-type inequality.

\begin{lemma}\label{lem:nash}
Let $\Psi$ be a $N$-function
and $\theta$ be an increasing function.
Assume that there exist a constant $\lambda >0$ such
that for any $x \geq 0$, $\theta(x/9) \geq \lambda \theta(x)$.
Let $\mu$ be a probability measure on $\dR^n$.

Then, the following are equivalent:\\
\noindent $(i)$
There exists a constant $C$ such that for
any smooth function $f: \dR^n \rightarrow \dR$ one has
$$
\var_\mu(f)
\theta \left( \frac12 \frac{\var_\mu(f)}{\NRM{f-\mu(f)}_\Psi^2} \right)
\leq  C \int |\nabla f|^2  d \mu .
$$
\noindent $(ii)$
There exists a constant $C'$ such that for
 any smooth function $f: \dR^n \rightarrow \dR$ one has
\begin{equation} \label{eq:mnash}
\var_\mu(f)
\theta \left( \frac12 \frac{\var_\mu(f)}{\NRM{f}_\Psi^2} \right)
\leq  C' \int |\nabla f|^2  d \mu .
\end{equation}
Furthermore, $(i) \Rightarrow (ii)$ with
$C' \leq C/\lambda$
and $(ii) \Rightarrow (i)$ with $C \leq C'$.
\end{lemma}

\begin{proof}
The implication $(ii)$ implies $(i)$ is obvious.

We will show that
$(i) \Rightarrow (ii)$. By \eqref{in:holderp},
for any function $f$,
$\int |f| d\mu \leq 2 \NRM{f}_\Psi \NRM{\ind}_{\Psi^*}
= \frac{2}{{\Psi^*}^{-1}(1)} \NRM{f}_\Psi$. It follows that
\begin{eqnarray*}
\NRM{f-\mu(f)}_\Psi
& \leq &
\NRM{f}_\Psi +  \NRM{\mu(|f|)}_\Psi
=
\NRM{f}_\Psi + \frac{1}{\Psi^{-1}(1)} \mu(|f|) \\
& \leq &
\left(1+\frac{2}{\Psi^{-1}(1){\Psi^*}^{-1}(1)}\right) \NRM{f}_\Psi
\leq 3 \NRM{f}_\Psi .
\end{eqnarray*}
In the last line we used the general bound
$x \leq \Psi^{-1}(x) {\Psi^*}^{-1}(x)$.
Since $\theta$ is increasing and $\theta (x/9) \geq \lambda \theta (x)$, it follows that
\begin{eqnarray*}
\var_\mu(f)
\theta \left( \frac12 \frac{\var_\mu(f)}{\NRM{f-\mu(f)}_\Psi^2} \right)
& \geq &
\var_\mu(f)
\theta \left( \frac1{18} \frac{\var_\mu(f)}{\NRM{f}_\Psi^2} \right) \\
& \geq &
\lambda \var_\mu(f)
 \theta \left( \frac1{2} \frac{\var_\mu(f)}{\NRM{f}_\Psi^2} \right) .
\end{eqnarray*}
Applying point $(ii)$ completes the proof.
\end{proof}

The second step is to link the Nash-type inequality in its simplified form
to an inequality between measure and capacity.

\begin{theorem}\label{th:nc}
Let $\Phi$ and $\Psi(x)=\frac{x^2}{\psi(|x|)}$
be two $N$-functions with $\psi$ increasing, satisfying $\psi(0)=0$
and $\lim_{+\infty}\psi =+\infty$.
let $\theta= {\Phi^*}^{-1} \circ \Psi \circ \psi^{-1}$.
Assume that\\
\noindent $(i)$ $x \mapsto \Psi \circ \psi^{-1}(x^2)$ is a Young function,\\
\noindent $(ii)$ there exists a constant
$\lambda >0$ such that for any $x \geq 0$,
$\theta(x/16) \geq \lambda \theta(x)$,\\
\noindent $(iii)$
there exists $\lambda' \ge 4 $ such that for all
$x\ge 2$ one has ${\Phi^*}^{-1}(\lambda' x)\le \lambda' {\Phi^*}^{-1}(x)/4$. \\
\noindent $(iv)$ the probability measure $\mu$ on $\dR^n$
satisfies for any smooth function $f: \dR^n \rightarrow \dR$
$$
\var_\mu(f)
\theta \left( \frac12 \frac{\var_\mu(f)}{\NRM{f}_\Psi^2} \right)
\leq  C \int |\nabla f|^2  d \mu ,
$$
for some constant $C$.

Then, for any Borel set $A$ such that
$\mu(A) < \frac 12$,
$$
\frac{1}{\Phi^{-1}\left( 1 / \mu(A) \right)} \leq
\frac{8 \lambda' C }{\lambda}
 \mathrm{Cap}_\mu(A) .
$$
\end{theorem}

\begin{proof}
Fix $A \subset \dR^n$ such that $\mu(A) < \frac 12$, and let
$g : \dR^n \rightarrow \dR$ such that
$g \geq \ind_A$ and $\mu(g=0) \geq \frac 12$.
Then for any $k \in \dZ$ we define $g_k=(g-2^k)_+\wedge 2^k$.
Let $H(x):=\Psi \circ \psi^{-1} (x^2)$. Note that
$\sqrt x H^{-1}(x) = \Psi^{-1}(x)$. Thus, by \eqref{in:holder} and
we have
$$
\NRM{g_k}_\Psi = \NRM{g_k \ind_{g_k \neq 0}}_\Psi
\leq 2 \NRM{g_k}_2 \NRM{\ind_{g_k \neq 0}}_H
= 2 \frac{\NRM{g_k}_2}{H^{-1}(1/\mu(g_k \neq 0))} .
$$
Note that $\mu(g_k = 0) = \mu(g \leq 2^k) \geq \mu(g=0) \geq \frac 12$.
Thus,
$$
\mu(g_k)^2 = \mu(g_k \ind_{\{g_k \neq 0\}})^2
\leq \mu(g_k^2) \mu(g_k \neq 0) \leq \frac 12 \mu(g_k^2)
$$
which in turn implies $\NRM{g_k}_2^2 \leq 2 \var_\mu (g_k)$.
This together with $\mu(g_k \neq 0) \leq \mu(g \geq 2^k)$ give
$$
\NRM{g_k}_\Psi^2 \leq 8
\frac{\var_\mu (g_k)}{[H^{-1}(1/\mu(g \geq 2^k))]^2} .
$$
Applying the Nash-type inequality to $g_k$ and the monotonicity of $\theta$, we get
\begin{eqnarray*}
\var_\mu (g_k) \theta
\left( \frac{1}{16} [H^{-1}(1/\mu(g \geq 2^k))]^2 \right)
& \leq &
\var_\mu (g_k)
\left( \frac 12 \frac{\var_\mu (g_k)}{\NRM{g_k}_\Psi^2} \right) \\
& \leq &
C \int |\nabla g|^2  d \mu .
\end{eqnarray*}
On the other hand,
\begin{eqnarray*}
\var_\mu (g_k)
& \geq &
\frac 12 \NRM{g_k}_2^2
\geq
2^{2k-1}
\NRM{\ind_{\{ g_k \geq 2^k\}}}_2^2
=
2^{2k-1} \mu(g_k \geq 2^k) \\
& = &
2^{2k-1} \mu(g \geq 2^{k+1}) .
\end{eqnarray*}
Let $\Omega_k = \{x : g(x) \geq 2^k \}$, $k \in \dZ$.
It follows from condition $(ii)$ on $\theta$
that for any $k \in \dZ$
\begin{eqnarray*}
\lambda 2^{2k-1} \mu(\Omega_{k+1})
\theta \left[ H^{-1}\left[\frac{1}{\mu(\Omega_k)} \right]^2 \right]
& \leq  &
2^{2k-1}
\mu(\Omega_{k+1})
\theta
\left[ \frac{1}{16} H^{-1}\left[\frac{1}{\mu(\Omega_k)} \right]^2 \right] \\
& \leq &
C \int |\nabla g|^2  d \mu .
\end{eqnarray*}
%
Now note that by definition of $H$ and $\theta$,
$\theta \left( H^{-1}(x)^2 \right)={\Phi^*}^{-1}(x)$.
Hence,
$$
\lambda 2^{2k-1}
\mu(\Omega_{k+1}) {\Phi^*}^{-1} \left(1/\mu(\Omega_k) \right)
\leq
C \int |\nabla g|^2  d \mu
\qquad \forall k \in \dZ .
$$
At this stage we may use \cite[Lemma 23]{barthe-cattiaux-roberto} we recall below with
$a_k =\mu(\Omega_k)$ and $F={\Phi^*}^{-1}$.
Since $\Phi^*$ is a Young function, the slope
function $x \mapsto \Phi^*(x)/x$ is non decreasing. This is equivalent to
say that $x \mapsto F(x)/x$ is non increasing. Thus the assumptions
of Lemma \ref{lem:ak} are satisfied, thanks to point $(iii)$. It follows that
$$
\lambda 2^{2k-1}
\mu(\Omega_{k}) {\Phi^*}^{-1} \left(1/\mu(\Omega_k) \right)
\leq
\lambda' C \int |\nabla g|^2  d \mu
\qquad \forall k \in \dZ .
$$
Furthermore,
by \eqref{eq:phiphi}, ${\Phi^*}^{-1}(x) \geq x/\Phi^{-1}(x)$. Hence,
$$
\lambda 2^{2k-1}
\frac{1}{{\Phi}^{-1} \left(1/\mu(\Omega_k) \right)}
\leq
\lambda' C \int |\nabla g|^2  d \mu
\qquad \forall k \in \dZ .
$$
Now take the largest $k$ such that $2^{2k} \leq1$. For that index,
$A \subset \{g \geq 2^k \} = \Omega_k$. By monotonicity it follows that
(using $1 \leq 2^{2(k+1})$)
$$
\frac{1}{\Phi^{-1}(1/\mu(A))}
\leq
2^{2(k+1)} \frac{1}{\Phi^{-1}(1/\mu(\Omega_{k}))}
\leq
\frac{8 \lambda' }{\lambda} C \int |\nabla g|^2  d \mu .
$$
The result follows by definition of the capacity.
\end{proof}

\begin{lemma}[\cite{barthe-cattiaux-roberto}] \label{lem:ak}
Let $ F: [2,+\infty) \to [0,+\infty)$ be a non-decreasing function
   such that $x\to F(x)/x$ is non increasing and there exists $\lambda'
   \ge 4 $ such that for all $x\ge 2$ one has $F(\lambda' x)\le \lambda'
   F(x)/4$.  Let $(a_k)_{k\in \dZ}$ be a non-increasing (double-sided)
   sequence of numbers in $[0,1/2]$. Assume that for all $k\in\dZ$ with
   $a_k >0$ one has
   $$
   2^{2k}a_{k+1}F\left(1/a_k\right) \le C,
   $$
   then for all
   $k\in \dZ$ with $a_k>0$ one has
   $$
   2^{2k}a_k F\left(1/a_k\right) \le \lambda' C.
   $$
\end{lemma}

We are now in position to give the
following reciprocal of Theorem
\ref{th:nash}.

\begin{corollary} \label{cor:nc}
Let $\Phi$ and $\Psi(x)=\frac{x^2}{\psi(|x|)}$
be two $N$-functions with $\psi$ increasing, satisfying $\psi(0)=0$
and $\lim_{+\infty}\psi =+\infty$.
let $\theta= {\Phi^*}^{-1} \circ \Psi \circ \psi^{-1}$.
Assume that\\
\noindent $(i)$ $x \mapsto \Psi \circ \psi^{-1}(x^2)$ is a Young function,\\
\noindent $(ii)$ there exists a constant
$\lambda >0$ such that for any $x \geq 0$,
$\theta(x/16) \geq \lambda \theta(x)$,\\
\noindent $(iii)$
there exists $\lambda' \ge 4 $ such that for all
$x\ge 2$ one has ${\Phi^*}^{-1}(\lambda' x)\le \lambda' {\Phi^*}^{-1}(x)/4$. \\
\noindent $(iv)$ the probability measure $\mu$ on $\dR^n$
satisfies for any smooth function $f: \dR^n \rightarrow \dR$
$$
\var_\mu(f)
\theta \left( \frac12 \frac{\var_\mu(f)}{\NRM{f-\mu(f)}_\Psi^2} \right)
\leq  C \int |\nabla f|^2  d \mu ,
$$
for some constant $C$.

Fix $k \in (0,+\infty)$ such that
for any $f$ with $f^2 \in \L_\Phi (\mu)$,
$\NRME{\mu(f)^2}{\Phi} \leq k \NRME{f^2}{\Phi}$.
Then, for any function $f$,
$$
\NRM{(f - \mu(f))^2}_\Phi
\leq \frac{64(1+k) \lambda'}{\lambda^2}
C \int |\nabla f|^2  d \mu .
$$
\end{corollary}

\begin{proof}
Apply Lemma \ref{lem:nash}, then Theorem \ref{th:nc}, and
finally Theorem \ref{th:cos} (together with Remark \ref{rem:thcos}).
\end{proof}



\section{Decay to equilibrium and Nash-type inequality}

Throughout this section we consider
a probability measure $d\mu=e^{-V(x)}dx$ on $\dR^n$
associated to a differentiable potential $V$ (or a limit of such measures).
Let $\GI=\Delta - \nabla V \cdot \nabla$
be a symmetric in $\L_2(\mu)$ diffusion generator
and $(\PT{t})_{t \geq 0}$ its associated semi-group.
In this setup we prove that Nash-type inequalities
are equivalent to the decay to equilibrium
of the semi-group in suitable Orlicz spaces associated
to the measure $\mu$.

\begin{theorem}\label{th:dnash}
Let $\Phi$ and $\Psi$
be two $N$-functions and $\theta$ an
increasing function.
Assume that the probability measure $\mu$ on $\dR^n$
satisfies, for any smooth function $f: \dR^n \rightarrow \dR$,
$$
\var_\mu(f)
\theta \left( \frac12 \frac{\var_\mu(f)}{\NRM{f-\mu(f)}_\Psi^2} \right)
\leq  C \int |\nabla f|^2  d \mu ,
$$
for some constant $C$.
Then, for any $t>0$,
$$
\var_\mu(\PT{t}f) \leq  m(t) \NRM{f-\mu(f)}_\Psi^2 ,
$$
where $m: \dR_+ \rightarrow \dR_+$ is the solution
of the differential equation
$$
m' = -\frac{2}{C_\Phi}
\theta \left( \frac{m}{2} \right)
$$
on $(0,\infty)$ such that $m(0)=+\infty$, or equivalently
$m$ satisfies for any $t\geq 0$,
$\displaystyle \int_{m(t)}^\infty \frac{1}{x\theta(x/2)} dx = \frac{2t}{C_\Phi}$.
\end{theorem}

\begin{proof}
Let $f$ be such that $\mu(f)=0$ and $\NRM{f}_\Psi=1$.
Set $u(t) = \var_\mu (\PT{t}f)$.
A differentiation gives
$$
u'(t) =
-2  \int |\nabla \PT{t}f|^2  d \mu
\leq
- \frac{2}{C_\Phi} u(t) \theta
\left( \frac 12 \frac{u(t)}{\NRM{\PT{t}f}_\Psi^2} \right) .
$$
Note that by convexity, $\NRM{\PT{t}f}_\Psi \leq \NRM{f}_\Psi = 1$.
Since $\theta$ is increasing we get
$$
u'(t) \leq - \frac{2}{C_\Phi}  u(t) \theta \left( \frac{u(t)}{2} \right) .
$$
By integration this gives
$$
\int_{u(t)}^{u(0)} \frac{dx}{x \theta(x/2)}
\geq \frac{2}{C_\Phi}t .
$$
Now, since
$$
\int_{m(t)}^\infty \frac{dx}{x \theta (x/2)}
= \frac{2}{C_\Phi}t ,
$$
we have that $u \leq m$ and the result
follows by homogeneity.
\end{proof}

Note that $m$ is not explicit in general. However
we can apply the Theorem to explicit examples.

\begin{example}
Let $\alpha \in [1,2]$,
$\beta = 2(1-\frac1\alpha) \in [0,1]$
and define  the probability measure on $\dR^n$:
$d\mu_\alpha^n(x) =
Z_\alpha^{-n} \exp \{-\sum_{i=1}^n |x_i|^\alpha \}dx$.
For any $\gamma \geq 1$ define
$\Phi_{\beta}^\gamma (x) = |x| \log(\gamma+|x|)^\beta$.
For $x \geq 0$ and $\delta \in (0,1)$, let also
$\psi(x)=\log(1+x)^\delta$ and $\Psi(x):=x^2/\psi(x)$.
We proved in Example \ref{ex:nash} that the following
Nash-type inequality holds: any $f$ satisfies
$$
\var_{\mu_\alpha^n}(f)^{1+\frac{\beta}{\delta}}
\leq C
\NRM{f-\mu_\alpha^n(f)}_{\Psi}^\frac{\beta}{\delta}
\int |\nabla f|^2  d \mu_\alpha^n .
$$
On the other hand, for $\theta=x^\frac\beta \delta$,
$$
\int_{m(t)}^\infty \frac{dx}{x \theta(x/2)}
=
\frac{2^\frac\beta \delta \delta}{\beta} \frac{1}{m(t)^\frac\beta \delta} .
$$
Hence, by Theorem \ref{th:dnash}
$$
\var_{\mu_\alpha^n}(\PT{t}f) \leq  2 \left( \frac{\delta C}{2\beta} \right)^\frac{\delta}{\beta}
\frac{1}{t^\frac{\delta}{\beta}}
\NRM{f-\mu(f)}_\Psi^2 .
$$
In other words, $\PT{t}$ is a continuous operator
from $\L_\Psi$ onto $\L_2$ with
$$
\NRM{\PT{t}}_{\L_\Psi \rightarrow \L_2} \leq
2 \left( \frac{\delta C}{2\beta} \right)^\frac{\delta}{\beta}
\frac{1}{t^\frac{\delta}{\beta}} .
$$
\end{example}

\begin{example}
As before, let $\alpha \in [1,2]$,
$\beta = 2(1-\frac1\alpha) \in [0,1]$ and
$d\mu_\alpha^n(x) =
Z_\alpha^{-n} \exp \{-\sum_{i=1}^n |x_i|^\alpha \}dx$ be a
probability measure on $\dR^n$.
For any $\gamma \geq 1$ define
$\Phi_{\beta}^\gamma (x) = |x| \log(\gamma+|x|)^\beta$.
Let also $\widetilde \psi(x)= e^{\log(1+x)^\delta}-1$ for
$\delta \in (0,1)$ and
$\widetilde \Psi(x)=x^2/\widetilde \psi(x)$.
We proved in Example \ref{ex:nash} that there exists a constant $C$ such that
for any $f$,
$$
\var_{\mu_\alpha^n}(f) \log \left( 1+ \frac 12
\frac{\var_{\mu_\alpha^n}(f)}
{\NRM{f-\mu_\alpha^n(f)}_{\Psi}^2} \right)^\frac{\beta}{\delta}
\leq C
\int |\nabla f|^2  d \mu_\alpha^n .
$$
For $\widetilde \theta (x): = \log(1+x)^\frac \beta \delta$, we define $m(t)$
as the unique solution of $2t/C = \int_{m(t)}^\infty dx/ [x \widetilde \theta (x/2)]$.
Now we deal with small values of $t$, small in such a way that $m(t) \geq 2$.
For such $t$'s we have
$$
\int_{m(t)}^\infty \frac{dx}{x \widetilde \theta(x/2)}
\leq
\int_{m(t)}^\infty \frac{dx}{x \log(\frac{x}{2})^\frac \beta \delta}
=
\frac{\delta}{\beta-\delta} \log\left(\frac{m(t)}{2}\right)^{\frac{\delta - \beta}{\delta}}
$$
provided that $\delta < \beta$ (if $\delta \geq \beta$ then $m$ is not defined!).
Hence,
$m(t) \leq \left( \frac{\delta C}{2(\beta - \delta)}
\frac{1}{t} \right)^{\frac{\delta}{\beta - \delta }}$.
Theorem \ref{th:dnash} implies that for small values of $t$,
$$
\var_{\mu_\alpha^n}(\PT{t}f) \leq  2
e^{C'/t^\frac{\delta}{\delta - \beta}}
\NRM{f-\mu(f)}_{\widetilde \Psi}^2
$$
with $C':=\left(\frac{\delta C}{2(\beta - \delta)} \right)^\frac{\delta}{\beta - \delta }$.
\end{example}

Now for completeness we prove a converse of Theorem \ref{th:dnash}.

\begin{theorem}\label{th:coul}
Let $\Psi$ be a $N$-function.
Assume that for any $t >0$,
$$
\var_\mu(\PT{t}f) \leq  m(t) \NRM{f-\mu(f)}_\Psi^2 .
$$
Then, for any smooth function
$f: \dR^n \rightarrow \dR$
$$
\var_\mu(f)
\widetilde{\theta}
\left( \frac12
\frac{\var_\mu(f)}{\NRM{f-\mu(f)}_\Psi^2} \right)
\leq
2 \int |\nabla f|^2  d \mu ,
$$
for $\widetilde{\theta}(x):=\sup_{t>0}\frac{1}{t}
\log \left( \frac{2x}{m(t)} \right)$.
\end{theorem}

\begin{proof}
We follow \cite[Proposition II.2]{coulhon96}.
Assume that $\mu(f)=0$ and that
$\NRM{f}_\Psi=1$. Let
$\int_0^\infty \lambda dE_\lambda$
be a spectral resolution of $-\GI$. Then
$\PT{t}=\int_0^\infty e^{-\lambda t} dE_\lambda$.
Since $\int_0^\infty
\frac{\mu(f \cdot dE_\lambda f)}{\mu(f^2)}=1$,
Jensen inequality yields
$$
\exp \left\{ \int_0^\infty (-2\lambda t)
\frac{\mu(f \cdot dE_\lambda f)}{\mu(f^2)} \right\}
\leq
\int_0^\infty e^{-2 \lambda t}
\frac{\mu(f \cdot dE_\lambda f)}{\mu(f^2)}.
$$
This exactly means that
$$
\exp \left\{-2t \frac{\mu(f \cdot (-\GI) f)}{\mu(f^2)}
\right\}
\leq
\frac{\NRME{\PT{t}f}{2}}{\mu(f^2)} .
$$
Now, using our assumption,
$\NRME{\PT{t}f}{2} \leq m(t)$.
Hence
$$
\frac{\mu(f^2)}{2t}
\log \left( \frac{\mu(f^2)}{m(t)} \right)
\leq
\mu(f \cdot (-\GI) f) = \int |\nabla f|^2 d\mu .
$$
The result follows by homogeneity and translation
invariance of the Dirichlet form.
\end{proof}

Next we recall a result, due to Grigor'yan,
which shows the link between
$m'(m^{-1})$ and $\widetilde{\theta}$.

We use the following definition, (\textit{cf.}
\cite{coulhon96}).
We say that a differentiable function
$m: (0,\infty) \to \dR^*_+$ satisfies condition $(D)$
if the derivative of its logarithm has polynomial growth,
{\it i.e.} $M(t)=-\log m(t)$ is such that
$$
M'(u) \geq \gamma M'(t),
\qquad \forall t>0, \quad \forall u \in [t,2t] ,
$$
for some $\gamma >0$ (for instance if $m$ behaves like
$t^{-d}$ or $e^{-ct^\delta}$, $0 \leq \delta \leq 1$
for $t$ large, it
satisfies condition $(D)$).

\begin{proposition}[\cite{grigoryan}]\label{prop:gre}
Let $m$ be a decreasing differentiable bijection
of $\dR^*_+$ satisfying condition $(D)$ for some
$\gamma>0$. Then, for all $x >0$,
$$
\widetilde{\theta}(x) =
\sup_{t>0} \frac{1}{t}
\log \left( \frac{2x}{m(t)} \right)
\geq
- \gamma \frac{m'\left(m^{-1} ( 2x) \right)}{x} .
$$
\end{proposition}

The above results imply the following equivalence of
the Nash-type inequality and the decay to equilibrium of the
semi-group.

\begin{theorem}\label{th:nasheq}
Let $m$ be a ${\cal C}^1$ decreasing bijection of
$\dR^*_+$ satisfying condition $(D)$ with $\gamma>0$.
Assume that $m'$ is an increasing bijection from
$\dR^*_+$ onto $\dR^*_-$.
Let
$\theta(|x|)= - m' \left( m^{-1}(2|x|)\right) / x$. Let
$\Phi$ and  $\Psi(x)=\frac{x^2}{\psi(|x|)}$
be two $N$-functions with $\psi$ increasing, satisfying $\psi(0)=0$
and $\lim_{+\infty}\psi =+\infty$. Assume that
$\theta= {\Phi^*}^{-1} \circ \Psi \circ \psi^{-1}$.

Then the following
are equivalent\\
$(i)$
for any $t >0$,
$$
\var_\mu(\PT{t}f) \leq  m(t) \NRM{f-\mu(f)}_\Phi^2 ,
$$
$(ii)$
for any smooth function $f: \dR^n \rightarrow \dR$
$$
\var_\mu(f)
\theta
\left( \frac12
\frac{\var_\mu(f)}{\NRM{f-\mu(f)}_\Psi^2} \right)
\leq
C_\Phi \int |\nabla f|^2  d \mu ,
$$
Moreover $(i) \Rightarrow (ii)$ with $C_\Phi=2/\gamma$
while $(ii) \Rightarrow (i)$ if $C_\Phi=2$.
\end{theorem}

\begin{proof}
To show that $(i)$ implies $(ii)$ it is enough to
apply Theorem \ref{th:coul} and proposition
\ref{prop:gre}.
The second part is a direct application of Theorem
\ref{th:dnash}.
\end{proof}





\begin{figure}
\psfrag{1}[]{$ \scriptstyle
\stackrel{\mathrm{Beckner-type } \; (T) \eqref{in:gb}:}{
\scriptstyle
\sup_{p \in (1,2)} \frac{\int f^2 d\mu - (\int f^p d\mu)^{2/p}}{T(2-p)}
\leq C_T \int |\nabla f|^2 d\mu}
$}
\psfrag{2}[][][1.5][90]{$\Longleftarrow \!\!=\!\!\! \Longrightarrow$ }
\psfrag{2b}[]{\scriptsize \cite{barthe-cattiaux-roberto}
(Theorem \ref{th:cgb})}
\psfrag{2c}[]{$
\genfrac{}{}{0pt}{2}{\Psi(x)= \frac{x}{T(1/\log(1+\frac{1}{x}))}}{}
$}
\psfrag{3}[]{$\scriptstyle \forall A \;\; \mathrm{ such \; that } \;\; \mu(A) < \frac 12,
\quad \mbox{Cap}_\mu(A) \ge C_\psi \psi(\mu(A))$}
\psfrag{4}[][][1.5][90]{$\Longleftarrow \!=\!=\!=\!=\!\!\!\!=\!= \!
\Longrightarrow$}
\psfrag{4b}{$\scriptstyle \Psi(x) = x$}
\psfrag{5}[][][1.5][90]{$ \Longleftarrow\!=\!=\!\!\!\!\! =$}
\psfrag{5b}{$
\begin{array}{l}
\mathrm{\scriptsize \cite[\scriptstyle Theorem 12]{barthe-cattiaux-roberto}} \\
\scriptstyle \Phi(x) =x \varphi(x) \; \mathrm{ and } \\
\scriptstyle \Psi(x)=x\varphi(\frac2x) 
\end{array}
$}
\psfrag{5c}{$\scriptstyle \stackrel{\mathrm{additive} \;
\Phi-\mathrm{Sobolev}\ ({\textbf{$\Phi$-S}})
:}
{\scriptstyle \int \! \Phi(f^2) d\mu -
\Phi( \int \!\! f^2 \! d\mu)\leq  c \! \int \! |\nabla  \!f|^2 d\mu
}$}
\psfrag{5d}[][][1.5][90]{$\Longleftarrow$}
\psfrag{5e}{\scriptsize Proposition \ref{prop:phi}}
\psfrag{5f}[][][1.5][90]{$\Longleftarrow
\!=\!=\!=\!=\!=\!=\!=\!=\!=\!=\!=\!=\!$}
\psfrag{5g}[][][1][90]{$\scriptstyle \mathrm {if} \; x \Phi''(x) \geq k$}
\psfrag{6}[][][1.5][90]{$\Longleftarrow\!=\!=\!=\!=\!=\!=\!\!\Longrightarrow$}
\psfrag{6b}{$
\begin{array}{l}
\scriptstyle \mathrm{Theorem \ref{th:cos}} \\
\scriptstyle \Psi(x)=\frac{1}{\Phi^{-1}(1/x)}
\end{array}
$}
\psfrag{6c}{$\scriptstyle \stackrel{\mathrm{Orlicz-Sobolev:}}{
\!\!\NRM{(f-\mu(f))^2}_\Phi \leq c \! \int \! |\nabla  \!f|^2 d\mu}
$}
\psfrag{6d}[][][1.5][90]{$\Longleftarrow$}
\psfrag{6e}{$\scriptstyle \mathrm{if} \; \forall f \;
\mu(|f|) \leq k  \NRM{f}_\Phi$}
\psfrag{7}[]{$\scriptstyle \stackrel{\mbox{\tiny 
{Poincar{\'e}}:}
}
{\genfrac{}{}{0pt}{2}{}{\var_\mu(f) \leq C_p \int |\nabla f|^2 d\mu}}$}
\psfrag{7b}[][][1.5][90]{$ \Longleftarrow\!=\!=\!=\!=\!=\!=$}
\psfrag{7c}[][][1][90]{\scriptsize Proposition \ref{prop:poincare}}
\psfrag{8}[][][1.5][270]{$\Longleftarrow \!\! = $}
\psfrag{8b}{\scriptsize Theorem \ref{th:nc}}
\psfrag{8c}{$\scriptstyle
\;\;\;\;\;\;\;\;\;\;\;  \stackrel{\mathrm{Nash-type} \; \eqref{eq:nash}:}
{\var_\mu(f) \theta \!
\bigl( \!  \frac{\var_\mu(f)}{2 \NRM{\widetilde f}_\Psi^2} \! \bigr)
\leq c \int \! |\nabla \! f|^2 d\mu}  $}
\psfrag{8d}[]{$\scriptstyle
\stackrel{\mathrm{modified \; Nash-type} \; \eqref{eq:mnash}:}
{\var_\mu(f) \theta \!
\bigl( \!  \frac{\var_\mu(f)}{2 \NRM{f}_\Psi^2} \! \bigr)
\leq c \int  |\nabla  f|^2 d\mu} $}
\psfrag{9}[][][1.5][90]{$\Longleftarrow \!\!\!\!\! \Longrightarrow$}
\psfrag{9b}{\scriptsize Lemma \ref{lem:nash}}
\psfrag{10}[][][1.5][180]{$\Longleftarrow$}
\psfrag{10b}{$
\begin{array}{l}
\scriptstyle \mathrm{Theorem \; \ref{th:nash}} \\
\end{array}
$}
\psfrag{11}[]{$\scriptstyle
\stackrel{\Phi-\mathrm{Entropy\; bound} \; ({\textbf{Ent$\Phi$-S}}):
}
{\int \!\Phi(f) d\mu - \Phi(\int\! f d\mu) \leq c \int \! |\nabla  f|^2
\Phi''(f) d\mu}
$}
\psfrag{11b}[][][1.5][270]{$ \Longleftarrow$}
\psfrag{11c}[]{\scriptsize \cite{chafai}}
\psfrag{12}[]{$\scriptstyle \mathrm{Ent}_\mu^\Phi( \PT{t}f) \leq e^{-mt}
\mathrm{Ent}_\mu^\Phi(f)$}
\psfrag{12b}[][][1.5][90]{$\Longleftarrow \!\!\!\!\! \!\!\!\Longrightarrow$}
\psfrag{12c}[]{\scriptsize \cite{chafai}}
\psfrag{13}[][][1.5]{$\Longleftarrow \!=$}
\psfrag{13b}[]{$\scriptstyle \mathrm{Take} \; p=1$}
\psfrag{14}[]{$\scriptstyle \stackrel{\mbox{\tiny{Poincar{\'e}}:}}
{\genfrac{}{}{0pt}{2}{}{\var_\mu(f) \leq C_p \int |\nabla f|^2 d\mu}}$}
\psfrag{15}[][][1.5][90]{$\Longleftarrow \!\!\!\!\! \Longrightarrow$}
\psfrag{15b}{$\scriptstyle \psi(x) =x$}
\psfrag{16}[][][1.5][90]{$\Longleftarrow \!\!\!\!\!\!\!\! \Longrightarrow$}
\psfrag{16b}{\scriptsize Theorem \ref{th:decay}}
\psfrag{16c}{\scriptsize \;\;\;\;\;\;\;\;\;\;
$\NRM{\PT{t} \widetilde f }_\Phi^2 \leq e^{-mt}
\NRM{ \widetilde f}_\Phi^2$}
\psfrag{17}[][][1.5][90]{$\Longleftarrow \!\! \!=\!=\!=\!=\!=\!=\!=$}
\psfrag{17b}[][][1][90]{$\scriptstyle \Phi'' \geq k>0$}
\psfrag{17c}{$\scriptstyle \;\;\;\;\;\;\;\;
\stackrel{\mathrm{modified \; Orlicz \; - \; Sobolev} \; \eqref{in:mos}:}
{\NRM{\widetilde f}_\phi^2 \leq c \int \! |\nabla  f|^2
\Phi''\bigl( \frac{ \widetilde f }{\NRM{\widetilde f}_\phi} \bigr)
d\mu}
$}
\psfrag{17d}[][][1.5][180]{$\Longleftarrow$}
\psfrag{17e}{\scriptsize Proposition \ref{prop:c2}}
\psfrag{18}[][][1.5][90]{$\Longleftarrow
 \!\!\!\!\!\!\!\! \Longrightarrow$}
\psfrag{18b}{\scriptsize Theorem \ref{th:nasheq} (and \ref{th:dnash}) }
\psfrag{18c}{$\scriptstyle \;\;\; \var_\mu( \PT{t} f ) \leq m(t)
\NRM{f-\mu(f)}_\phi^2$}
\epsfig{file=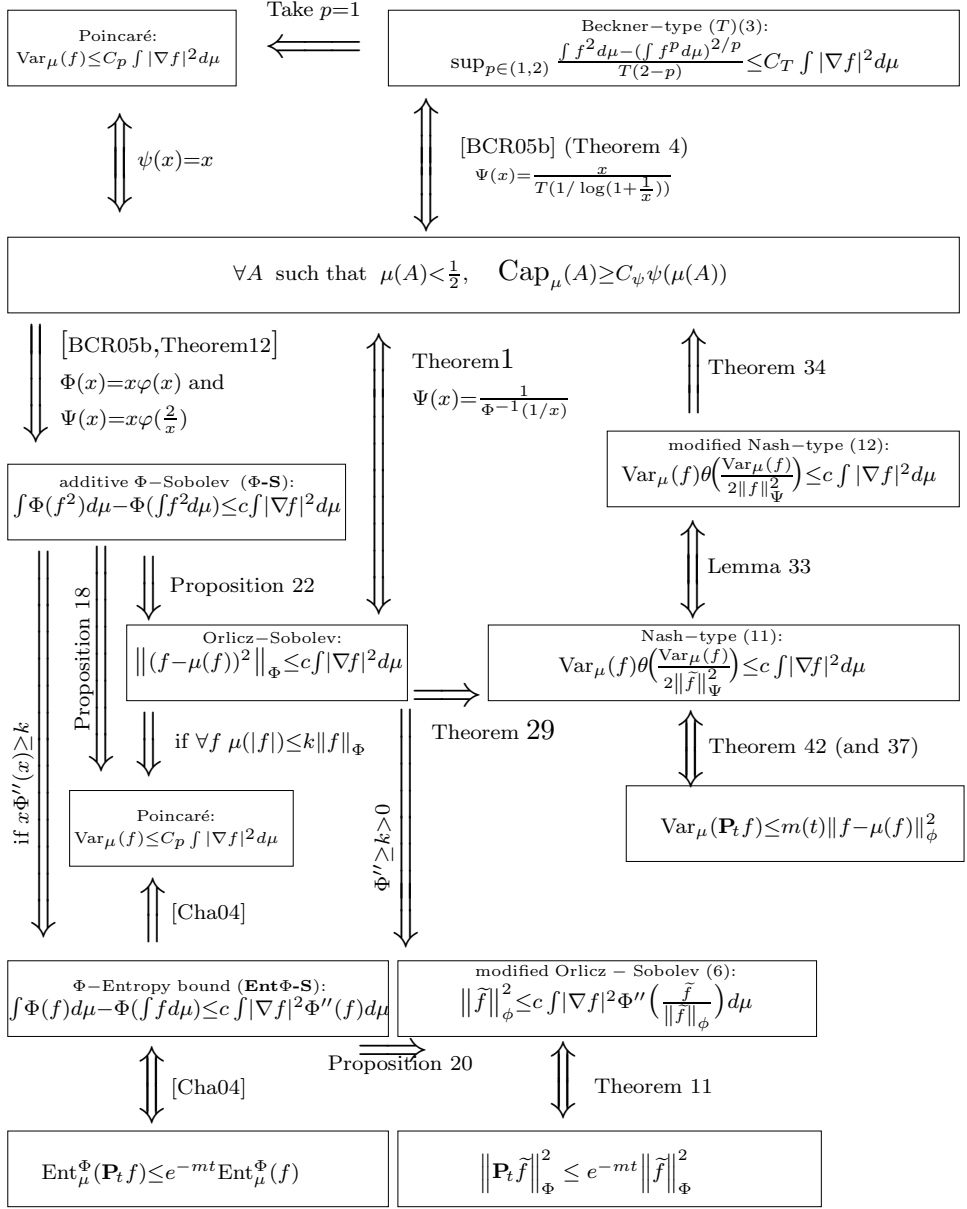,height=16cm,width=12.6cm,angle=0}
\caption{Implications Network. Here $\widetilde f:=f-\mu(f)$.}
\end{figure}

\newpage

\section{Appendix : Young functions and Orlicz spaces}
\label{sec:orlicz}

In this section we collect some results
on Orlicz spaces.
We refer the reader to \cite{rao-ren} for
demonstrations and complements.

\begin{definition}[Young function]
A function $\Phi : \dR \rightarrow [0,\infty]$
is a \emph{Young function}
if it is convex, even, such that $\Phi(0)=0$, and
$\lim_{x \rightarrow +\infty} \Phi(x)=+\infty$.
\end{definition}

The Legendre transform $\Phi^*$ of $\Phi$ defined by
$$
\Phi^*(y)= \sup_{x \geq 0} \{ x|y| - \Phi(x) \}
$$
is a lower semi-continuous Young function.
It is called the \emph{complementary function} or
\emph{conjugate} of $\Phi$.

Among the Young functions, we will consider
those continuous
with finite values such that
$\Phi(x)/x \rightarrow \infty$
as $x \rightarrow \infty$
(for stability reasons w.r.t. duality).
When additionally $\Phi(x)=0 \Leftrightarrow x=0$ and
$\Phi'(0_+)=0$, $\Phi$ is called a  \emph{$N$-function}.

For any lower semi-continuous Young function
$\Phi$ (in particular
if $\Phi$ has finite values), the conjugate of
$\Phi^*$ is $\Phi$.
The pair $(\Phi,\Psi)$ is said to be a
\emph{complementary pair} if
$\Psi = \Phi^*$ (or equivalently $\Phi = \Psi^*$).
 When $\Phi(1) + \Phi^*(1)=1$,
the pair $(\Phi, \Phi^*)$ is said to be
\emph{normalized}. The conjugate of an
$N$-function is an $N$-function.

Let $\Phi$ be an $N$-function. Then,
for any $a>0$,
\begin{equation} \label{eq:phiphi}
a < \Phi^{-1}(a) (\Phi^*)^{-1}(a) \leq 2a .
\end{equation}
\smallskip
The simplest example of N-function is
$\Phi(x)=\frac{|x|^p}{p}$,
$p>1$, in which
case, $\Phi^*(x)=\frac{|x|^q}{q}$, with $1/p+1/q=1$.
The function $\Phi(x)=|x|^\alpha \ln(1+|x|)^\beta$
is also a Young function for
$\alpha \geq 1$ and $\beta \geq 0$
and an $N$-function when
$\alpha > 1$ or $\beta >0$.

Now let $(\cal X,\mu)$ be a measurable space,
and $\Phi$ a Young function.
The space
$$
\L_\Phi(\mu)=\{f: {\cal X} \rightarrow \dR
\mbox{ measurable};
\exists \alpha > 0,
\int_{\cal X} \Phi(\alpha f) < + \infty \}
$$
is called the {\it Orlicz space} associated to $\Phi$.
When $\Phi(x)=|x|^p$, then $\L_\Phi(\mu)$ is the standard
Lebesgue space  $\L_{p}(\mu)$.

There exist two equivalent norms which give to $\L_\Phi(\mu)$ a structure of
Banach space. Namely, \textit{Luxemburg norm}
$$
\NRME{f}{\Phi} = \inf\{ \lambda >0; \int_{\cal X}
\Phi \left( \frac{f}{\lambda} \right) d\mu \leq 1 \}
$$
and \textit{the Orlicz norm}
$$
N_\Phi(f)=\sup \{\int_{\cal X} |fg|d\mu; \int_{\cal X}
\Phi^*(g) d\mu \leq 1 \} \;.
$$
Note that we invert the notation with respect
to \cite{rao-ren}.
We will use the notation ${\cal G}_\Phi$,
or more simply $\cal G$
when no confusion, the set
${\cal G}_\Phi = \{|g|:
\int_{\cal X} \Phi^*(g) d\mu \leq 1\}$. Note
in particular that ${\cal G}_\Phi$ is a
space of non negative functions.
Moreover
\begin{equation} \label{in:eqn}
\NRME{f}{\Phi}
\leq
N_\Phi(f)
\leq
2 \NRME{f}{\Phi} .
\end{equation}
By definition of the norm and the previous result, it is easy to see that
for any measurable subset $A$ of $\cal X$,
\begin{equation} \label{equiv}
\NRM{\ind_A}_{\Phi}
=
\frac{1}{\Phi^{-1} \left( \frac{1}{\mu(A)} \right)} .
\end{equation}

Then, the following result generalizes H{\"o}lder inequality.
Let $\Phi_1$, $\Phi_2$ and $\Phi_3$ be three Young functions satisfying
for all $x \geq 0$, $\Phi_1^{-1}(x) \Phi_2^{-1}(x) \leq \Phi_3^{-1}(x)$.
Then, for any $(f,g) \in \L_{\Phi_1}(\mu) \times \L_{\Phi_2}(\mu)$,
\begin{equation}\label{in:holder}
\NRME{fg}{\Phi_3}
\leq
2 \NRME{f}{\Phi_1}\NRME{g}{\Phi_2} \;.
\end{equation}
In particular, when $\Phi_3(x)=|x|$, we get
$\int_{\cal X} |fg| d\mu
\leq 2 \NRME{f}{\Phi_1}\NRME{g}{\Phi_2}$.
In the case of complementary pairs of
Young functions, we have the
following more precise result, see
\cite[Proposition 1 in section 3]{rao-ren}:
\begin{equation}\label{in:holderp}
\int_{\cal X} |fg|  d\mu \leq
2\NRME{f}{\Phi}\NRME{g}{\Phi^*} .
\end{equation}

Finally, for any constant $c >0$, it is easy to see that for any function $f$,
\begin{equation}\label{in:cst}
c \NRME{f}{\Phi(\cdot/c)} = \NRME{f}{\Phi} .
\end{equation}

\subsection*{Comparison of norms}

Let us notice that any Young function $\Phi$ satisfies
$\left| x \right| = O \left( \Phi(x) \right)$
as $x$ goes to $\infty$. It leads to the following lemma.

\begin{lemma} \label{lemma_injection_L1}
Any Orlicz space may
be continuously embedded in $\L_1$. More precisely,
let $D$ and $\tau$ in $(0,\infty)$ such that
$
\left| x \right| \leq \tau \,  \Phi(x)
$
for any $\left| x \right| \geq D$. Then, for any
$f \in \L_\Phi$,
\begin{equation}
\label{norm_injection}
\NRM{f}_1 \leq (D+\tau) \, \NRM{f}_\Phi.
\end{equation}
Consequently, if $\Phi$ and $\Psi$ are two Young functions satisfying,
for some constants $A,B \geq 0$,
$\Phi(x) \leq A  |x| + B  \Psi(x)$,
then
\begin{equation}
\label{norm_comparison}
\NRM{f}_\Phi \leq \max
\left(1 , A \NRM{\text{\emph{Id}}}_{\L_\Psi
\rightarrow \L_1} + B  \right)  \NRM{f}_\Psi.
\end{equation}
\end{lemma}

\begin{remark}
\label{remark_norm_mu_f2}
When $\Phi(x)/x \rightarrow \infty$ as
$x\rightarrow \infty$,
we may choose
$\tau = 1$ or any other positive constant. We get in
particular the estimate
\begin{equation}
\label{bound_norm_mu_f2}
\NRM{\mu (f)^2}_\Phi \leq \left(D + 1 \right)
\NRM{\ind}_\Phi  \NRM{f^2}_\Phi,
\end{equation}
where $D$ is such that $|x| \leq \Phi(x)$ for any $|x| \geq D$.
\end{remark}

\begin{proof}[Proof of lemma \ref{lemma_injection_L1}]
Let $f \in \L_\Phi(\mu)$. We may assume by homogeneity  that
$\NRM{f}_\Phi =1$.
Then $ \int \Phi(f) \, d \mu = 1$ and so
\begin{eqnarray*}
\int \left| f \right|  d \mu
& = &
\int_{\{\left| f \right| \leq D\}}
| f| d \mu + \int_{\{|f| \geq D\}}
|f| d \mu\\
& \leq&
D \mu \left(| f| \leq D \right) + \tau
\int_{\{|f| \geq D\}} \Phi(f) d \mu \leq D + \tau.
\end{eqnarray*}
As for bound \eqref{norm_injection}, assume now
that $\NRM{f}_\Psi =1$ and hence
$\int \Psi(f) \, d \mu = 1$ as well.
For any $\lambda \geq 1$,
\begin{eqnarray*}
\int \Phi \left(\frac{f}{\lambda} \right) d \mu
&\leq&
\frac{A}{\lambda} \NRM{f}_1 +
B \int \Psi \left(\frac{f}{\lambda} \right) d \mu \\
& \leq &
\frac{A}{\lambda}
\NRM{\text{Id}}_{\L_\Psi \rightarrow \L_1}
\NRM{f}_\Psi + \frac{B}{\lambda}
\int \Psi(f) d \mu \leq 1
\end{eqnarray*}
provided
$
\lambda
\geq A \NRM{\text{Id}}_{\L_\Psi \rightarrow \L_1} + B$.
Note that for the second inequality we used
convexity of $\Psi$.
\end{proof}



\newcommand{\etalchar}[1]{$^{#1}$}
\def\cprime{$'$}

\end{document}